\documentclass[12pt]{amsart}
\usepackage[colorlinks=true,citecolor=blue,linkcolor=blue]{hyperref}
\usepackage{amssymb,amscd,amsmath,graphicx,mathtools}
\usepackage{enumerate}
\usepackage{fullpage}
\usepackage{xcolor}
\usepackage{amsthm}
\usepackage{float}
\usepackage{pgfplots}
\usepackage{dirtytalk}
\usepackage{listings}
\usepackage{longtable}
\usepackage{mathrsfs}
\usepackage{bbm}
\usepackage{comment}
\usepackage{graphicx}

\newtheorem{thm}{Theorem}[section]

\newtheorem{prop}[thm]{Proposition}
\newtheorem{lem}[thm]{Lemma}

\newtheorem{quest}[thm]{Question}

\newtheorem{stp}[thm]{Setup}

\theoremstyle{definition}
\newtheorem{defn}[thm]{Definition}

\newtheorem{rem}[thm]{Remark}

\newtheorem{exmp}[thm]{Example}

\theoremstyle{remark}

\definecolor{energy}{RGB}{114,0,172}
\definecolor{freq}{RGB}{45,177,93}
\definecolor{spin}{RGB}{251,0,29}
\definecolor{signal}{RGB}{203,23,206}
\definecolor{circle}{RGB}{217,86,16}
\definecolor{average}{RGB}{203,23,206}

\colorlet{shadecolor}{gray!20}
\pgfplotsset{compat=1.9}

\usepgflibrary{fpu}

\newcommand{\mm}{{\mathfrak m}}
\newcommand{\ZZ}{{\mathbb Z}}
\newcommand{\NN}{{\mathbb N}}
\newcommand{\QQ}{{\mathbb Q}}
\newcommand{\RR}{{\mathbb R}}

\begin{document}
	
\title{Limits of length functions of multi $p$-families of ideals}

\author{Th\'ai Th\`anh Nguy$\tilde{\text{\^E}}$n}
\address{McMaster University, Department of Mathematics and Statistics,
	1280 Main Street West, Hamilton, Ontario, Canada \\
	and University of Education, Hue University, 34 Le Loi St., Hue, Viet Nam}
\email{nguyt161@mcmaster.ca  tnguyen11@tulane.edu}

\author{Vinh Anh Ph{\d{A}}m}
\address{Tulane University \\ Department of Mathematics \\
	6823 St. Charles Ave. \\ New Orleans, LA 70118, USA}
\email{vpham1@tulane.edu}


\begin{abstract}
We show the asymptotic relationship between the limit of the normalized length function of a multi-$p-$family of ideals and that of its shifted family under linear growth conditions in a local domain of characteristic $p$. Examples of multi-$p-$families of ideals including products of Frobenius powers of different ideals. We apply our results to obtain a generalized version of a formula due to Wantanabe-Yoshida for certain $p-$families using results from Verma, and to provide an instance of the existence of a mixed multiplicity version of multi-$p$-families of ideals.
\end{abstract}

\maketitle

%%%%%%%%%%%%%%%%%%%%%%%%%%%%%%%%%%

%%%%%%%%%%%%%%%%%%%%%%%%%%%%%%%

\section{Introduction} 
\label{sec.intro}

It is well-known that for a $\mm$-primary ideal $I$ in a commutative Noetherian local ring $R$, the Hilbert-Samuel function $n\to \lambda_R(R/I^{n})$, where $\lambda(M)$ is the length function of the module $M$, coincides with a polynomial for large enough values of $n$. However, if one considers the function $q\to \lambda_R(R/I^{[q]})$, where $q$ is a power of $p$ and $I^{[q]}=\langle x^q\ |\ x\in I \rangle$, the $q$-th Frobenius power of $I$, in the ring $R$ of characteristic $p>0$, its behavior is much more complicated. Kunz's paper \cite{Kunz76} can be seen as one of the first to study the function $q\to \lambda_R(R/I^{[q]})$, namely, the Hilbert-Kunz function of $I$. Remarkably, Monsky \cite{Mon83} showed that the limit
\begin{equation*}
    e_{HK}(I,R):=\lim\limits_{e\to \infty} \dfrac{\lambda_R(R/I^{[q]})}{q^d}
\end{equation*}
exists for any $\mathfrak{m}$-primary ideal $I$. The above limit is called the Hilbert-Kunz multiplicity of $I$. Much current research concerns the Hilbert-Kunz function and the Hilbert-Kunz multiplicity, for example, see \cite{DJ18,T18,S19,GKV21} and the references therein. Besides, the definition of Frobenius's power of an ideal in prime characteristic can be extended to arbitrary non-negative real exponents in \cite{HTW20}. This plays an important role in studying the singularities of algebraic varieties.

The Newton-Okounkov body theory offers a powerful method for studying the asymptotic behaviors and invariants in algebraic geometry and commutative algebra. Originated in the seminal work of Okounkov \cite{Ok96}, this construction was systematically developed and utilized in a study of multiplicity in \cite{LM2009,KK2012}. In \cite{KK14}, Kaveh and Khovanskii used similar construction of convex bodies associated to graded families of ideals in regular local rings and other nice local rings to study their multiplicity, extending earlier results in \cite{RS78,Mustata02} with the same spirit. A collection $\{ I_n\}_{n\in \NN}$ of ideals is called a graded family if $I_pI_q\subseteq I_{p+q}$ for all $p,q\in \NN$. Later, Cutkosky successfully extended the method Newton-Okounkov bodies to study the numerical limits of length functions and multiplicities of graded families of ideals in a more general setting of local rings (\cite{C13,C14,C15}). For other uses of the Newton-Okounkov body in studying asymptotic behaviors and properties in commutative algebra, see \cite{Mayes14a,HN23} and the references therein. In prime characteristic, an analogous concept namely the $p$-bodies was introduced by D. J. Hern\'andez and J. Jeffries in \cite{DJ18} and used to study length functions as well as the Hilbert-Kunz multiplicity and other related invariants of \emph{$p$-families of ideals} in local rings. A collection of ideals $\mathcal{I}=\{I_q\}_{q=1}^{\infty}$ in $R$ indexed by powers of $p$ is called a $p$-family if $I_q^{[p]}\subseteq I_{pq}$ for all $q$. Recently, much attention has been drawn to $p$-families in terms of the similarities and differences between them and graded families, especially about multiplicities. Many results showed the existence of the mixed multiplicities of a multi-graded family of ideals. However, the existence of a positive characteristic version of the mixed multiplicities of a multi-$p$- family is still unresolved. 

Cutkosky-Sarkar-Srinivasan \cite{CSS19} initiated the study of mixed multiplicities of graded families of ideals for the case of $\mathfrak{m}-$primary filtration. Recently, under mild assumptions, namely the \emph{linear growth conditions} (see below for details), the notion of mixed multiplicities has been extended to graded families of $\mathfrak{m}$-primary ideals by the work of Cid-Ruiz and Monta\~no \cite{YJ22}. They also observed some important results such as a \say{Volume = Multiplicity formula} for mixed multiplicities of graded families. In their study, one of the main tools is to show the asymptotic relationship between the limit of the normalized length function of a multi-graded family of ideals and that of its shifted family under mild conditions. More specifically, let $\{J(1)\}_{n\in \mathbb{N}},\ldots, \{J(r)\}_{n\in \mathbb{N}}$ be graded families of non-zero ideals, and let $\{I(1)\}_{n\in \mathbb{N}},\ldots, \{I(s)\}_{n\in \mathbb{N}}$ be $\mathfrak{m}$-primary graded families of ideals, for $\mathbf{n}=\{n_1,\ldots,n_r\}\in \mathbb{N}^r$ and $\mathbf{m}=\{m_1,\ldots,m_s\}\in \mathbb{N}^s$ the two following limits exist and are equal 
\begin{equation*}
	\lim\limits_{p\to \infty}\lim\limits_{m\to \infty}\dfrac{\lambda(\mathbf{J}(p)^{m\mathbf{n}}/\mathbf{I}(p)^{m\mathbf{m}}\mathbf{J}(p)^{m\mathbf{n}})}{p^{d}m^d}=\lim\limits_{m\to \infty}\dfrac{\lambda(\mathbf{J}_{m\mathbf{n}}/\mathbf{I}_{m\mathbf{m}}\mathbf{J}_{m\mathbf{n}})}{m^d},
\end{equation*}
where
\begin{align*}
   &\mathbf{J}(p)^{m\mathbf{n}}:= J(1)_p^{mn_1}\cdots J(r)_p^{mn_r}, \mathbf{I}(p)^{m\mathbf{m}}:= I(1)_p^{mm_1}\cdots I(s)_p^{mm_s},\\
   &\mathbf{J}_{m\mathbf{n}}:= J(1)_{mn_1}\cdots J(r)_{mn_r}, \text{\space and \space} \mathbf{I}_{m\mathbf{n}}:= I(1)_{mm_1}\cdots I(r)_{mm_s}.
\end{align*}

In our study, we proved that the result above still holds for the case of $p$-families of ideals and derive some interesting applications. We will work with the following setup.

Let $(R,\mathfrak{m}, \mathbbm{k})$ be a Noetherian ring of dimension $d$ and characteristic $p>0$ with perfect residue field $\mathbbm{k}$, $\hat{R}$ denote the $\mathfrak{m}$-adic completion of $R$. Let $\{J(1)_q\}_{q=1}^{\infty},\ldots,\{J(r)_q\}_{q=1}^{\infty}$ be $p$-families of non-zero ideals, and let $\{I(1)_q\}_{q=1}^{\infty},\ldots,\{I(s)_q\}_{q=1}^{\infty}$ be $\mathfrak{m}$-primary $p$-families of ideals. For $\underline{n}=(n_1,\ldots,n_r)\in \mathbb{N}^r,\underline{m}=(m_1,\ldots,m_s)\in \mathbb{N}^s$, and $q=p^e, e\in \NN$, we use the following notation:
	\begin{equation*}
		\begin{split}
			\mathbf{J}_{q{\underline{n}}}&:=J(1)_{qp^{n_1}}\cdots J(r)_{qp^{n_r}}, \mathbf{I}_{q\underline{m}}:=I(1)_{qp^{m_1}}\cdots I(s)_{qp^{m_s}},\\
			\mathbf{J}(p^b)^{q\underline{n}}&:=J(1)_{p^b}^{[qp^{n_1}]}\cdots J(r)_{p^b}^{[qp^{n_r}]}, \mathbf{I}(p^b)^{q\underline{m}}:=I(1)_{p^b}^{[qp^{m_1}]}\cdots I(s)_{p^b}^{[qp^{m_s}]}.
		\end{split}
	\end{equation*}
We also define the pair of $p$-families
\begin{equation*}
	\begin{split}
	(\mathscr{J}_{\underline{n}},\mathscr{H}_{\underline{m},\underline{n}})&:=\left(\{\mathbf{J}_{q\underline{n}}\}_{q=1}^{\infty}, \{\mathbf{I}_{q\underline{m}}\mathbf{J}_{q\underline{n}}\}_{q=1}^{\infty}\right) \text{\space and\space}\\
	(\mathscr{J}(p^b)_{\underline{n}},\mathscr{H}(p^b)_{\underline{m},\underline{n}})&:=\left(\{\mathbf{J}(p^b)^{q\underline{n}}\}_{q=1}^{\infty}, \{\mathbf{I}(p^b)^{q\underline{m}}\mathbf{J}(p^b)^{q\underline{n}}\}_{q=1}^{\infty}\right) \text{\space and\space} \text{\space for every\space} b\in \mathbb{N}.
	\end{split}
\end{equation*}

\noindent We say that the pair of $p$-families $(\mathcal{J},\mathcal{I})$ has \emph{linear growth} if there exists $c=c(\mathcal{J},\mathcal{I})\in \mathbb{N}$ such that 
	\begin{equation*}
		\mathfrak{m}^{cq}\cap J_q=\mathfrak{m}^{cq}\cap I_q\text{\space for every power\space} q \text{\space of\space} p.  
	\end{equation*} 

In our first main result, we assume the following linear growth conditions on the pair the above pair families: each $(\mathscr{J}_{\underline{n}},\mathscr{H}_{\underline{m},\underline{n}})$ has linear growth, and that if $c_{\underline{m},\underline{n}}:=c(\mathscr{J}_{\underline{n}},\mathscr{H}_{\underline{m},\underline{n}})$, then $c(\mathscr{J}(p^b)_{\underline{n}},\mathscr{H}(p^b)_{\underline{m},\underline{n}})=c_{\underline{m},\underline{n}}.p^b$ for every $b\in \mathbb{N}$. It is worth pointing out that the linear growth condition is a natural condition, and automatically holds, for example, when $\mathcal{J}(1), \mathcal{J}(2), \ldots ,\mathcal{J}(r)$ are the ring $R$. We now state our first main result that provides us an important technical tool to study length functions and multiplicities of multi-$p$-families.

\begin{thm}[Theorems \ref{maintool}, \ref{thm.maintool}]
    \label{mainThmSec3}
    Let $(R,\mathfrak{m},\mathbbm{k})$ be a Noetherian local ring of dimension $d$ and characteristic $p>0$ such that $\dim (\mathrm{N}(\hat{R}))<d$; here $\mathrm{N}(\hat{R})$ denotes the nilradical of the $\mathfrak{m}$-adic completion $\hat{R}$.  We have that the following limits exist and are equal 
	\begin{equation*}
		\lim\limits_{b\to \infty}\lim\limits_{q\to \infty}\dfrac{\lambda(\mathbf{J}(p^b)^{q\underline{n}}/\mathbf{I}(p^b)^{q\underline{m}}\mathbf{J}(p^b)^{q\underline{n}})}{p^{bd}q^d}=\lim\limits_{q\to \infty}\dfrac{\lambda(\mathbf{J}_{q\underline{n}}/\mathbf{I}_{q\underline{m}}\mathbf{J}_{q\underline{n}})}{q^d}.
	\end{equation*}
\end{thm}

Our first application of Theorem \ref{mainThmSec3} is a generalization of results given in \cite{WY01} concerning the length function of multi-$p$-families, their multiplicities, mixed multiplicities, and Hilbert-Kunz multiplicities in a $2$-dimensional Cohen-Macaulay ring.

\begin{thm}
\label{mainThmSec41}
    Let $(R,\mathfrak{m})$ be a $2$-dimensional Cohen-Macaulay local ring over characteristic $p>0$ such that $\dim (\mathrm{N}(\hat{R}))<d$, where $\mathrm{N}(\hat{R})$ denotes the nilradical of the $\mathfrak{m}$-adic completion $\hat{R}$. Let $\mathcal{I}(1)=\{I(1)_n\}_{n=1}^{\infty},\ldots,\mathcal{I}(s)=\{I(s)_n\}_{n=1}^{\infty}$ be $\mathfrak{m}$-primary $p$-families of ideals such that 
 \begin{equation}
  I(i)_q^{[p]}\subseteq I(i)_q^p\subseteq I(i)_{qp}\text{\space for every \space} i=1,\ldots,s.   
 \end{equation}
 For every $b\in \mathbb{N}$, assume that $r(I(i)_{p^b}|I(j)_{p^b})=0$ for $1\leq i\leq j\leq s$. Then for any $\underline{m}=(m_1,\ldots,m_s)\in \mathbb{N}^s$, we have 
 {\footnotesize
    \begin{equation*}
        \lim\limits_{q\to \infty}\dfrac{\lambda(R/\mathbf{I}_{q\underline{m}})}{q^2}=\sum\limits_{i=1}^s\left(\lim\limits_{b\to \infty}\dfrac{e(I(i)_{p^b})}{p^{2b}}\binom{p^{m_i}}{2}+\lim\limits_{b\to \infty}\dfrac{e_{HK}(I(i)_{p^b})}{p^{2b}}p^{m_i}\right)+\sum\limits_{1\leq i\leq j\leq s}\lim\limits_{b\to \infty}\dfrac{e(I(i)_{p^b}|I(j)_{p^b})}{p^{2b}}p^{m_i}p^{m_j}.
    \end{equation*}}
\end{thm}

In particular, if $R$ is an analytically unramified, Cohen-Macaulay local rings of dimension $2$ and the families of ideals satisfy the condition in the above theorem, then the result holds. \par
\vspace{0.5em}

Our next application of Theorem \ref{mainThmSec3} is one instance where one can define by showing the existence of a mixed multiplicity version of multi-$p-$families of ideals.

\begin{thm}[Theorems \ref{thm.minsumclosure},\ref{thm.mixedMultCoeff}]
    \label{mainThmSec42}
    Let $(R,\mathfrak{m},\mathbbm{k})$ be a Noetherian local domain of dimension $d$ and characteristic $p>0$ with the fraction field $\mathbb{F}$. Consider a collection of ideals $\mathcal{M}$ in $R$ such that there exists a valuation $\vartheta$ so that $\vartheta(I\cdot J)=\vartheta(I)+\vartheta(J)$ for any ideals $I,J \in \mathcal{M}$. Let $\mathcal{I}(1) = \{I(1)_q\}_{q=1}^{\infty},\ldots,\mathcal{I}(1) = \{I(s)_q\}_{q=1}^{\infty}$ be $p$-families
    of ideals in such collection $\mathcal{M}$.

    Suppose that there is a non-degenerate linear transformation $T$ such that $T(\mathcal{C})\subseteq \RR^d_{\ge 0}$. Furthermore, for every $b \in \mathbb{N}$, suppose that the closure of the $p$-body  $${\Delta(S,\vartheta(I(i)_{p^b})}=\bigcup\limits_{q=1}^{\infty}\left(\dfrac{1}{q}\vartheta(I(i)^{[q]}_{p^b})+\mathrm{Cone}(S)\right)$$
    is a cobounded $C$-convex region. Then the limit
   \begin{equation*}
       \lim\limits_{q\to \infty}\dfrac{\lambda(R/\mathbf{I}(p^b)^{q\underline{n}})}{q^d}
   \end{equation*}
   coincides with a homogeneous polynomial in $p^{n_1},\ldots,p^{n_s}$ of total degree equal to $d$. 

   As a consequence, for each $b\in \NN$, set
    \[
    \mathcal{P}_{\mathbf{I}(p^b)} (p^{n_1},\ldots ,p^{n_s}) = \lim\limits_{q\to \infty}\dfrac{\lambda(R/\mathbf{I}(p^b)^{q\underline{n}})}{q^d}.
    \]
    Then there exists a homogeneous polynomial of total degree $d$ with real coefficients, denoted $\mathcal{P}_{\mathbf{I}}$, such that 
    \[
    \mathcal{P}_{\mathbf{I}}(p^{n_1},\ldots ,p^{n_s}) = \lim\limits_{b\to \infty} \mathcal{P}_{\mathbf{I}(p^b)} (p^{n_1},\ldots ,p^{n_s}) = \lim\limits_{q\to \infty}\dfrac{\lambda(R/\mathbf{I}_{q\underline{n}})}{q^d},
    \]
    for all $(n_1,\ldots,n_s)\in \NN$.
\end{thm}

%%%%%%%%%%%%%%%%%%%%%%%%%%%%%%%%%%%%%%%%
\par
\vspace{1em}
\noindent
{\bf Acknowledgments.}
The authors thank Huy T\`ai H\`a for 
his comments and suggestions.
The first author is partially supported by the NAFOSTED (Vietnam) under the grant number 101.04-2023.07. This is part of the second author's PhD thesis.
%%%%%%%%%%%%%%%%%%%%%%%%%%%%%%%%%%%%%%%%%%%

\section{Preliminaries}
\label{sec.prelim}
In this section, we recall terminology and results that will be used often in the paper. For other unexplained terminology from algebra and convex geometry, we refer the interested readers to the following texts \cite{DJ18,Zie1995}.

\subsection{Convex Geometry}
\label{subsec.convgeo}

A convex cone is a subset of $\mathbb{R}^d$ that is closed under taking an $\mathbb{R}$-linear combination of points with non-negative coefficients. Denote $\mathrm{Cone}(U)\subseteq \mathbb{R}^d$ the convex cone that is the closure of the set of all linear combinations $\sum_{i}\lambda_i u_i$ with $u_i\in U$ and $\lambda_i\in \mathbb{R}_{\geq 0}$. A cone in $\mathbb{R}^d$ has a non-empty interior if and only if the real vector space it generates has dimension $d$. We call such a cone full-dimensional. 

A cone is \emph{pointed} if it is closed and if there exists a vector $\mathbf{a}\in \mathbb{R}^d$ such that  $\langle \mathbf{u},\mathbf{a} \rangle >0$ for all $\mathbf{u}\in C\setminus \{0\}$, where $\langle - , - \rangle$ is the usual inner product in $\mathbb{R}^d$.

If $C$ is a pointed cone and $\alpha$ is a non-negative real number, we define:
\begin{equation*}
	H=H_{\alpha}:=\{\mathbf{u}\in \mathbb{R}^d \ | \ \langle \mathbf{u},\mathbf{a} \rangle <\alpha\}
\end{equation*} 
a \emph{truncating half-space} for $C$.

%%%%%%%%%%%
%%%%%%%%%%%%%%%%%%%%%%%%%%%%%%%%%%%
\subsection{Semigroups}
\label{subsec.val}

A semigroup is a subset of $\mathbb{Z}^d$ that contains $0$ and is closed under addition. A semigroup $S$ is finitely generated if there exists a finite subset $S_0\subseteq S$ such that every element of $S$ can be written as an $\mathbb{N}$-linear combination of elements of $S_0$. A subset $T\subseteq S$ of a semigroup $S$ is called an \emph{ideal} of $S$ if $T+S \subseteq T$. A semigroup $S$ is called \emph{pointed} if, whenever $\mathbf{a}\in S$ and $-\mathbf{a}\in S$, we must have $\mathbf{a} = \mathbf{0}$. A semigroup $S$ is called \emph{standard} if $S-S=\mathbb{Z}^d$, and the full dimensional cone generated by $S$ in $\mathbb{R}^d$ is pointed.

%%%%%%%%%%%%%
%%%%%%%%%%%%%%%%%%%%%%%%%%%%%%%%%%%%%
\subsection{$p$-systems and $p$-bodies}
\label{subsec.psysbod}

\begin{defn}\cite[Definition 4.3]{DJ18}
\label{def.p-system}
A collection of subsets $T_{\bullet}=\{T_{q}\}_{q=1}^{\infty}$ of a semigroup $S$ indexed by $q=p^e$ ($p$ is a prime) is called a $p$-system if it satisfies the following conditions 
	\begin{enumerate}
		\item $T_q$ is an ideal of $S$ for all $q$,
		\item $pT_q\subseteq T_{pq}$ for all $q$.
	\end{enumerate}
\end{defn}

\begin{defn}\cite[Definition 4.4]{DJ18}
	The $p$-body associated to a given $p$-system of ideals $T_{\bullet}$ of a semigroup $S$ in $\mathbb{Z}^d$ is defined by 
	\begin{equation*}
		\Delta(S,T_{\bullet})=\bigcup_{q=1}^{\infty}\dfrac{1}{q}T_q+\mathrm{Cone}(S).
	\end{equation*}
\end{defn}

We recall the following two useful results used extensively in our paper. These results can be considered the $p$-body versions of similar results in the Newton-Okounkov body theory \cite{KK2012, LM2009}. 

\begin{thm}\cite[Theorem 4.10]{DJ18}\label{limit=volumetheorem}
	For a standard semigroup $S$ in $\mathbb{Z}^d$, a $p$-system $T_{\bullet}$ in $S$, and a truncating halfspace $H$ for $\mathrm{Cone}(S)$, we have:
	\begin{equation*}
		\lim\limits_{q\to \infty}\dfrac{\#(T_q\cap qH)}{q^d}=\mathrm{Vol}_{\mathbb{R}^d}(\Delta(S,T_{\bullet})\cap H).
	\end{equation*}
\end{thm}

\begin{thm}\cite[Theorem 3.18]{SD22}\label{Fujitatypeapproximationtheorem}
	For a standard semigroup $S$ in $\mathbb{Z}^d$, a $p$-system $T_{\bullet}$ in $S$, and a truncating halfspace $H$ for $\mathrm{Cone}(S)$. Then for any $\varepsilon>0$, there exists $q_0$ such that if $q\geq q_0$ the following inequality holds 
	\begin{eqnarray*}
			\lim\limits_{e\to \infty}\dfrac{\#((p^eT_q+S)\cap p^eqH)}{p^{ed}q^d}\geq\mathrm{Vol}_{\mathbb{R}^d}(\Delta(S,T_{\bullet})\cap H)-\varepsilon.
	\end{eqnarray*}
\end{thm}

%%%%%%%%%%%%%
%%%%%%%%%%%%%%%%%%%%%%%%%%%%%%%%%%%%%
\subsection{Valuations and OK-valuations}\label{section_about_OKvaluation}

Let $\mathbb{F}$ be a field and $\mathbb{F}^{\times}=\mathbb{F}\setminus \{0\}$. Fix $\mathbf{a}\in \mathbb{R}^d$, and fix an embedding $\mathbb{Z}^d \hookrightarrow \mathbb{R}$ given by $\mathbf{v}\mapsto \langle \mathbf{a},\mathbf{v} \rangle$. An $\mathbf{a}$-valuation on a field $\mathbb{F}$ with the value group $\mathbb{Z}^d$ is a surjective group homomorphism $\vartheta:\mathbb{F}^{\times}\to \mathbb{Z}^d$ with the property that 
\begin{align*}
	\vartheta (xy)&=\vartheta(x)+\vartheta(y),\\
	\vartheta(x+y)&\geq_{\mathbf{a}}\min\{\vartheta(x),\vartheta(y)\}.
\end{align*} 
Here, by $\textbf{u} \ge_{\textbf{a}} \textbf{v}$, we mean $\langle \mathbf{a},\mathbf{u} \rangle \ge \langle \mathbf{a},\mathbf{v} \rangle$. For $N\subseteq \mathbb{F}$, let 
is $\vartheta(N):=\vartheta(N^{\times})$ the image of $N$ under $\vartheta$, where $N^{\times}=N\setminus\{0\}$. Given a point $\mathbf{u}\in \mathbb{Z}^d$, we define:
\begin{align*}
	\mathbb{F}_{\geq \mathbf{u}}:=\{x\in \mathbb{F}|\vartheta(x)\geq_{\mathbf{a}}\mathbf{u}\}\cup \{0\} \text{\space and\space} \mathbb{F}_{> \mathbf{u}}:=\{x\in \mathbb{F}|\vartheta(x)>_{\mathbf{a}}\mathbf{u}\}\cup\{0\}.
\end{align*}
The valuation ring $(V_{\vartheta},\mathfrak{m}_{\vartheta},\mathbbm{k}_{\vartheta})$ of $\mathbb{F}$ is associated to $\vartheta$, that is, $V_{\vartheta}=\mathbb{F}_{\geq \mathbf{0}}$, and $\mathfrak{m}_{\vartheta}=\mathbb{F}_{>\mathbf{0}}$. A local domain $(R,\mathfrak{m},\mathbbm{k})$ is said to be dominated by $(V_{\vartheta},\mathfrak{m}_{\vartheta},\mathbbm{k}_{\vartheta})$  if $(R,\mathfrak{m},\mathbbm{k})$ is a local subring of $(V_{\vartheta},\mathfrak{m}_{\vartheta},\mathbbm{k}_{\vartheta})$, i.e., $\mathfrak{m}\subseteq\mathfrak{m}_{\vartheta}$. 

Following the work of Cutkosky \cite{C13, C14}, Hern\'andez and Jeffries in \cite{DJ18} defined a distinguished class of valuations, namely, the OK-valuations, that extend the notion of “good valuations” with “1-dimensional leaves” that were defined and used by Kaveh and Khovanskii in their work on Newton-Okounkov body \cite{KK2012}.

\begin{defn}\cite[Definition 3.1]{DJ18}\label{OKrelativedef}
	Let $(R,\mathfrak{m},\mathbbm{k})$ be a local domain of dimension $d$ with fraction field $\mathbb{F}$, and fix a $\mathbb{Z}$-linear embedding of $\mathbb{Z}^d \hookrightarrow \mathbb{R}$. If a valuation
		$\vartheta:\mathbb{F}^{\times}\to \mathbb{Z}^d$
on $\mathbb{F}$ with the value group $\mathbb{Z}^d$ and local ring $(V_{\vartheta},\mathfrak{m}_{\vartheta},\mathbbm{k}_{\vartheta})$ satisfies the following conditions:
\begin{enumerate}
	\item $R$ is strongly dominated by $V_{\vartheta}$,
	\item the resulting extension of the residue fields $\mathbbm{k}\hookrightarrow\mathbbm{k}_{\vartheta}$ is finite, and 
	\item there exists a point $\mathbf{v}\in \mathbb{Z}^d$ such that for every $a \in \ZZ_{\ge 0}$, we have
		$R\cap \mathbb{F}_{\geq a\mathbf{v}}\subseteq\mathfrak{m}^a$,
\end{enumerate} 
then we say that $(V_{\vartheta},\mathfrak{m}_{\vartheta},\mathbbm{k}_{\vartheta})$ is OK relative to $R$. A local domain $R$ of dimensional $d$ is said to be OK if there exists a valuation on its fraction field with value group $\mathbb{Z}^d$ that is OK relative to $R$.	
\end{defn}

\begin{stp}\cite[Setup 5.8]{DJ18}\label{setuplocaldomainwithOKvaluation}
For the rest of this section, we fix a $d$-dimensional local domain $(R,\mathfrak{m},\mathbbm{k})$ of characteristic $p>0$ with the perfect residue field $\mathbbm{k}$ and the fraction field $\mathbb{F}$, a $\mathbb{Z}$-linear embedding $\mathbb{Z}^d \hookrightarrow \mathbb{R}$ induced by $\mathbf{a} \in \mathbb{R}^d$, and a valuation $\vartheta:\mathbb{F}^{\times} \twoheadrightarrow \mathbb{Z}^d$ that is OK relative to $R$. We use $S$ to denote the semigroup $\vartheta(R)$ in $\mathbb{Z}^d$, and $C$ to denote the closed cone in $\mathbb{R}^d$ generated by $S$. 
\end{stp}

Note that $\vartheta(R)$ is a standard semigroup. 

\begin{rem}
	By assumption, $R$ is a local subring of $V$. Therefore, if $M$ is an $R$-submodule of $\mathbb{F}$ and $\mathbf{u}\in \mathbb{Z}^d$, then the quotient
		$\dfrac{M\cap \mathbb{F}_{\geq \mathbf{u}}}{M\cap \mathbb{F}_{>\mathbf{u}}}$
has a $\mathbbm{k}$-vector space structure induced by the $R$-module structure on $M$. By definition, this space is nonzero if and only if there exists an element $m\in M$ with $\vartheta(m)=\mathbf{u}$. 
\end{rem}

\begin{defn}\cite[Definition 3.9]{DJ18}
	Let $M$ be a $R$-submodule of $\mathbb{F}$ and $1\leq h\leq [\mathbbm{k}_{\vartheta}:\mathbbm{k}]$, define 
	\begin{equation*}
		\vartheta^{(h)}(M)=\left\{\mathbf{u}\in \mathbb{Z}^d \ |\ \dim_{\mathbbm{k}}\left(\dfrac{M\cap \mathbb{F}_{\geq \mathbf{u}}}{M\cap \mathbb{F}_{>\mathbf{u}}}\right)\geq h\right\}.
	\end{equation*}
\end{defn}

\begin{rem}
	Let $M$ be a $R$-submodule of $\mathbb{F}$. For all $1\leq h\leq [\mathbbm{k}_{\vartheta}:\mathbbm{k}]$, $\vartheta^{(h)}(M)$ is an ideal of $S$. In fact, for any $g\in R^{\times}$, with $\mathbf{v}=\vartheta(g)$, and any $\mathbf{u}\in \mathbb{Z}^d$, the map \begin{equation*}
		\dfrac{M\cap \mathbb{F}_{\geq \mathbf{u}}}{M\cap \mathbb{F}_{>\mathbf{u}}}\to \dfrac{M\cap \mathbb{F}_{\geq \mathbf{u}+\mathbf{v}}}{M\cap \mathbb{F}_{>\mathbf{u}+\mathbf{v}}}
	\end{equation*}
defined by $[m]\mapsto[gm]$ is a $\mathbbm{k}$-linear injection. Therefore 
	$\vartheta^{(h)}(M)+S\subseteq\vartheta^{(h)}(M)$. 

\end{rem}

The following lemma is very useful to compute the length function of $M$ via the number of integral points in its associated semigroup.

\begin{lem}\cite[Lemma 3.11]{DJ18}\label{lengthcomputedbyimagesofMunderv(h)}
	If $M$ is a $R$-submodule of $\mathbb{F}$ and $\mathbf{v}\in \mathbb{Z}^d$, then 
	\begin{equation*}
		\lambda_D(M/M\cap \mathbb{F}_{\geq \mathbf{v}})=\sum\limits_{h=1}^{[\mathbbm{k}_{\vartheta}:\mathbbm{k}]}\#\left(\vartheta^{(h)}(M)\cap H\right),
	\end{equation*}
where $H$ is the halfspace $\{\mathbf{u}\in \mathbb{R}^d:\left \langle \mathbf{u},\mathbf{a} \right \rangle< \left \langle \mathbf{v},\mathbf{a} \right \rangle\}.$
\end{lem}

%%%%%%%%%%%%
%%%%%%%%%%%%%%%%%%%%%%%%%%%%%%%%%%%%%%
\subsection{$p$-families of ideals}
\label{subsec.pfam}

\begin{defn}\cite[Definition 5.1]{DJ18}
A sequence of ideals $\mathcal{I}=\{I_q\}_{q=1}^{\infty}$ in a ring $R$ indexed by powers of $p$ is called a $p$-family if $I_q^{[p]}\subseteq I_{pq}$ for all $q$, where $J^{[p]} = \langle x^{p} \ | \ x\in J \rangle$ denotes the $p$-th Frobenius power of the ideal $J$.	
\end{defn}

$p$-families of ideals are ubiquitous in characteristic-$p$ commutative algebra. Examples include the family of $p^e$-th Frobenius powers of an ideal, the termwise product, sum, or intersection of an arbitrary collection of $p$-families defines a $p$-family as well as termwise expansion, contraction via a map, or saturation with respect to some ideal of $p$-families (see \cite[Examples 5.6, 5.7]{DJ18}).

\begin{rem}
\label{rem.mprimarypropertyofpfamily}
	Given a $p$-family of ideals $\mathcal{I}$ in a local ring $(R,\mathfrak{m})$, the ideals in this family are $\mathfrak{m}$-primary if and only if $I_1$ is $\mathfrak{m}$-primary. Indeed, if $\mathfrak{m}^a$ is contained in $I_1$ for some positive integer $a$, and $I_1$ is generated by $b$ elements, then, by the pigeon-hole principle,
	\begin{equation*}
		\mathfrak{m}^{abq}\subseteq I_1^{bq}\subseteq I_1^{[q]}\subseteq I_q.
	\end{equation*}
	In other words, if $c=ab$, then $\mathfrak{m}^{cq}\subseteq I_q$ for all $q$ a power of $p$.
\end{rem}

Let $(R,\mathfrak{m},\mathbbm{k})$ be a $d$-dimensional local OK-domain with an OK-valuation $\vartheta$. We observe that if $x\in R$, then $p\cdot\vartheta(x)=\vartheta(x^p)$. Thus, if we choose any $p$-family $\mathcal{I}$ of ideals, then 
\begin{enumerate}
	\item $\vartheta(I_q)$ is an ideal of the semigroup $\vartheta(R)=S$, and 
	\item $p\vartheta(I_q)\subseteq\vartheta(I_q^{[p]})\subseteq \vartheta(I_{pq})$.
\end{enumerate}
Therefore $\{\vartheta(I_q)\}_{q=1}^{\infty}$ is a $p$-system of ideals in $S$.\par
\vspace{0.5em}

Note that $\vartheta^{(h)}(I_q)$ is an ideal of $S$ for every $q$, but the collection $\vartheta^{(h)}(\mathcal{I})=\{\vartheta^{(h)}(I_q)\}_{q=1}^{\infty}$ might not be a $p$-system unless we assume that $\mathbbm{k}$ is perfect, see \cite[Remark 5.9]{DJ18}. However, one still can uniformly approximate $\vartheta^{(h)}(\mathcal{I})$ by $\vartheta(\mathcal{I})$.

\begin{thm}\cite[Corollary 5.10]{DJ18}\label{limit=volumeforp-familiesthm}
	Let $R$ be a $d$-dimensional local OK domain with OK valuation $\vartheta$. For a $p$-family of ideals $\mathcal{I}$ in $D$ we have 
	\begin{equation*}
		\lim\limits_{q\to \infty}\dfrac{\#(\vartheta^{(h)}(I_q)\cap qH)}{q^d}=\mathrm{Vol}_{\mathbb{R}^d}(\Delta(S,\vartheta(\mathcal{I}))\cap H),
	\end{equation*}
where $S=\vartheta(R)$, $C=\mathrm{Cone}(S)$ and $H$ is any truncating half-space of $C$.
\end{thm}

%\textcolor{red}{It is worth to have the above theorem only for the case when $\vartheta^{(h)}(I_q)$ is not $p$-system because we already have theorem \ref{limit=volumetheorem}.}

We sometimes use the notation $[\vartheta^{(h)}(\mathcal{I})]_q$ for $\vartheta^{(h)}(I_q)$. We also have the following lemma whose proof can be deduced from \cite[Lemma 5.21]{DJ18}. 

\begin{lem}\label{lem.length.over.R.to.sum.length.over.Ri}
	Let $(R,\mathfrak{m},\mathbbm{k})$ be a $d$-dimensional reduced local ring of positive characteristic, and $\mathcal{I}, \mathcal{J}$ be sequences of ideals of $R$ indexed by the powers of $p$ such that $\mathfrak{m}^{cq}\subseteq I_q\subseteq J_q$ for some positive integer $c$ and for all $q=p^e, e\in \mathbb{N}$. Let $P_1,\ldots, P_n$ be the minimal primes of $R$ and $R_i=R/P_i$, then there exists $\gamma>0$ such that for all $q$,
	\begin{equation*}
		\left|\sum\limits_{i=1}^n\lambda_{R_i}(J_qR_i/I_qR_i)-\lambda_R(J_q/I_q)\right|\leq \gamma.q^{d-1}.
	\end{equation*}
\end{lem}

\begin{comment}
   We have the following exact sequences \textcolor{blue}{may not even need to mention these exact sequences}
    \begin{equation*}
        0\to J_q/I_q \to R/I_q\to R/J_q\to 0
    \end{equation*}
    and 
    \begin{equation*}
        0\to J_qR_i/I_qR_i \to R_i/I_qR_i\to R_i/J_qR_i\to 0
    \end{equation*}
    for every $i=1,\ldots,n$. 
\end{comment}

\begin{proof}
     It follows from \cite[Lemma 5.21]{DJ18} that there exist $\delta_1,\delta_2>0$ such that for all $q$,
    \begin{align*}
        &\left|\sum\limits_{i=1}^n\lambda_{R_i}(J_qR_i/I_qR_i)-\lambda_R(J_q/I_q)\right|\\
        =& \left|\sum\limits_{i=1}^n(\lambda_{R_i}(R_i/I_qR_i)-\lambda_{R_i}(R_i/J_qR_i))-(\lambda_R(R_q/I_q)-\lambda_R(R_q/J_q))\right|\\
        =& \left|\left(\sum\limits_{i=1}^n\lambda_{R_i}(R_i/I_qR_i)-\lambda_R(R_q/I_q)\right)-\left(\sum\limits_{i=1}^n\lambda_{R_i}(R_i/J_qR_i)-\lambda_R(R_q/J_q)\right)\right|\\
        \leq& \left|\sum\limits_{i=1}^n\lambda_{R_i}(R_i/I_qR_i)-\lambda_R(R_q/I_q)\right|+\left|\sum\limits_{i=1}^n\lambda_{R_i}(R_i/J_qR_i)-\lambda_R(R_q/J_q)\right|\\
        \leq & \delta_1.q^{d-1}+\delta_2.q^{d-1} = \gamma.q^{d-1},  \text{ where } \gamma:=\delta_1+\delta_2.
    \end{align*}
   
\end{proof}

%%%%%%%%%%%%%%
%%%%%%%%%%%%%%%%%%%%%%%%%%%%%%%%%%%%%
\section{Length Functions of $p-$ Families with Linear Growth Condition}
\label{sec.lingrowth}

The purpose of this section is to provide the main technical results for our applications in the next section. We show the existence of certain limits of length functions and relate the limits to the volume of $p$-families. Our methods utilize the standard treatment using truncation families and shifted families, which is analogous to that for graded families in \cite{C14,CSS19,YJ22}. We first fix the following setup.

\begin{stp}\label{pfamiliessetup}
	Let $(R,\mathfrak{m},\mathbbm{k})$ be a $d$-dimensional complete local domain of characteristic $p>0$. Let $\mathcal{J}=\{J_q\}_{q=1}^{\infty}$ and $\mathcal{I}=\{I_q\}_{q=1}^{\infty}$ be $p$-families of non-zero ideals, such that $J_q\supseteq I_q$ for every $q=p^e, e\in \mathbb{N}$. For every $a\in \mathbb{Z}_{\ge 0}$, define $\mathcal{J}_{p^a}=\{J_{p^a,q}\}_{q=1}^{\infty}$ to be the $p^a$-th truncated $p$-family, with the element
	\begin{equation*}
		J_{p^a,p^e}=\begin{cases}
			J_{p^e} &\text{\space if \space} e\leq a\\
			\sum\limits_{\tiny\begin{matrix} 
					i,j\geq0\\
					i+j=e
				\end{matrix}}J_{p^a,p^i}^{[p^{j}]}&\text{\space if \space} e>a. 
		\end{cases}
	\end{equation*}
Likewise, define $\mathcal{I}_{p^a}:=\{I_{p^a,q}\}_{q=1}^{\infty}$ given by the above rule.  
\end{stp}

From the definition and $p-$family property of $\mathcal{J}$, we have $J_{p^a,p^e} \subseteq J_{p^e}$ for all $e\in \mathbb{Z}_{\ge 0}$. Moreover, the truncated families possess the following special property.

\begin{defn}
	A $p$-family $\mathcal{I}=\{I_q\}_{q=1}^{\infty}$ is said to be \emph{of finite type} if there exists a power of $p$, says $q'$, such that $I_{q'}^{[p^e]}=I_{q'p^e}$ for all $e\in  \mathbb{Z}_{\ge 0}$.
\end{defn}

\begin{lem}
\label{lem.Noethtruncation}
    For every $a\in \mathbb{Z}_{\ge 0}$, the family $\mathcal{J}_{p^a}=\{J_{p^a,q}\}_{q=1}^{\infty}$ is a $p$-family of finite type.
\end{lem}

\begin{proof}
    For a fix $a$, by induction on $e \ge a+1$, we have 
    \[
    J_{p^a,p^e} = J_1^{[p^e]} + J_p^{[p^{e-1}]} + \cdots + J_{p^a}^{[p^{e-a}]}.
    \]
   In particular, for $e\ge a+1$, $J_{p^a,p^e}^{[p^b]}=J_{p^a,p^ep^b}$ for all $b\in \mathbb{Z}_{\ge 0}$. Thus, $\mathcal{J}_{p^a}$ is a $p$-family of finite type. 
\end{proof}

The following definition includes a condition on $p$-families that is needed in our main results. This is a natural condition in our context and was used in the context of graded families, see \cite[Theorem 6.1]{C14} and \cite{YJ22}.

\begin{defn}
	Assume Setup \ref{pfamiliessetup}. We say that the pair of $p$-families $(\mathcal{J},\mathcal{I})$ has linear growth if there exists $c=c(\mathcal{J},\mathcal{I})\in \mathbb{N}$ such that 
	\begin{equation*}
		\mathfrak{m}^{cq}\cap J_q=\mathfrak{m}^{cq}\cap I_q\text{\space for every power\space} q \text{\space of\space} p.  
	\end{equation*}  
\end{defn}

We note that the above condition holds when $\mathcal{J}$ and $\mathcal{I}$ are $\mm$-primary, see Remark \ref{rem.mprimarypropertyofpfamily}.\par
\vspace{0.5em}

Our first main result provides an approximation for the volume of the $p$-body of a $p$-family by that of its truncated families. This result can be thought of as the $p$-body version to \cite[Proposition 4.3]{CSS19} and \cite[Lemma 2.5]{YJ22}.

\begin{thm}	
\label{thm.volume=volumeoftruncation}
Adopt the context of Setup \ref{setuplocaldomainwithOKvaluation}, and fix a $p$-family of ideals $\mathcal{J}$ in $R$. Then we have 
	\begin{align*}
		&\lim\limits_{a\to \infty}\mathrm{Vol}_{\mathbb{R}^d}(\Delta(S,\vartheta^{(h)}(\mathcal{J}_{p^a}))\cap H)=\mathrm{Vol}_{\mathbb{R}^d}(\Delta(S,\vartheta^{(h)}(\mathcal{J}))\cap H).
		%,\text{\space and\space}\\
		%&\lim\limits_{a\to \infty}\mathrm{Vol}_{\mathbb{R}^d}(\Delta(S,\vartheta(\mathcal{I}_{p^a})))=\mathrm{Vol}_{\mathbb{R}^d}(\Delta(S,\vartheta(\mathcal{I})))
	\end{align*}
 for all $1\leq h\leq [\mathbbm{k}_{\vartheta}:\mathbbm{k}]$ and any truncating half-sapce $H$.
\end{thm}

\begin{proof}

Since $\mathbbm{k}$ is perfect, then $\vartheta^{(h)}(\mathcal{J})$ is a $p$-system of ideals in $S$ %\textcolor{red}{Do we really need $\vartheta^{(h)}(\mathcal{J})$ to be a p-system? Corrolary 5.10 in \cite{DJ18} does not need that to be true (do not need k to be perfect) We can still use the argument in Corrolary 4.17 of \cite{DJ18} to show a same version of Theorem \ref{Fujitatypeapproximationtheorem} without assuming $\vartheta^{(h)}(\mathcal{J})$ be a p-system, it that true?}. 
It follows that for every $e,a\in \mathbb{Z}_{>0}$ and $1\leq h\leq [\mathbbm{k}_{\vartheta}:\mathbbm{k}]$, we have $$p^e\vartheta^{(h)}(J_{p^a})+S\subseteq\vartheta^{(h)}(J_{p^a}^{[p^e]})+S\subseteq \vartheta^{(h)}(J_{p^a}^{[p^e]})\subseteq \vartheta^{(h)}(J_{p^a,p^ap^e}).$$ %and $p^e\vartheta^{(h)}(I_{p^a})+S\subseteq \vartheta^{(h)}(I_{p^a,p^ap^e})$. 
Thus, by Theorem \ref{limit=volumeforp-familiesthm} and Theorem \ref{Fujitatypeapproximationtheorem}, for a fixed $\varepsilon\in \mathbb{R}_{>0}$, there exists $a_0\in \mathbb{N}$ such that if $a\geq a_0$ and for any truncating halfspace $H$, we have 
\begin{equation}\label{sosanh1}\begin{split}
			\mathrm{Vol}_{\mathbb{R}^d}(\Delta(S,\vartheta^{(h)}(\mathcal{J}))\cap H)&\geq \mathrm{Vol}_{\mathbb{R}^d}(\Delta(S,\vartheta^{(h)}(\mathcal{J}_{p^a}))\cap H)\\
		&=\lim\limits_{e\to \infty}\dfrac{\#(\vartheta^{(h)}(J_{p^a,p^ap^e})\cap p^ap^eH)}{p^{ad}p^{ed}}\\
		&\geq \lim\limits_{e\to \infty}\dfrac{\#((p^e\vartheta^{(h)}(J_{p^a})+S)\cap p^ap^eH)}{p^{ad}p^{ed}}\\
		&\geq \mathrm{Vol}_{\mathbb{R}^d}(\Delta(S,\vartheta^{(h)}(\mathcal{J}))\cap H)-\varepsilon.
	\end{split}
\end{equation}  
The statement now follows as $\varepsilon$ is arbitrary.
\end{proof}

Before proceeding to our next estimation results, we fix the following setup of multi-$p$-families including the shifted families and the truncated families.

\begin{stp}\label{multipfamiliessetup}
	We adopt Setup \ref{pfamiliessetup} and further assume that the residue field $\mathbbm{k}$ is perfect. Let $\mathcal{J}(1)=\{J(1)_q\}_{q=1}^{\infty},\ldots,\mathcal{J}(r)=\{J(r)_q\}_{q=1}^{\infty}$ be $p$-families of non-zero ideals, and let $\mathcal{I}(1)=\{I(1)_q\}_{q=1}^{\infty},\ldots,\mathcal{I}(s)=\{I(s)_q\}_{q=1}^{\infty}$ be $p$-families of $\mathfrak{m}$-primary ideals. For $\underline{n}=(n_1,\ldots,n_r)\in \mathbb{N}^r,\underline{m}=(m_1,\ldots,m_s)\in \mathbb{N}^s$, $a,b\in \mathbb{N}$, and $q=p^e, e\in \NN$, we use the following notation:
 
%	\begin{align*}\label{notationformultifamiliessetup}
%		\begin{split}
%			\mathbf{J}_{q\underline{n}}&:=J(1)_{qp^{n_1}}\cdots J(r)_{qp^{n_r}}, \mathbf{I}_{q\underline{m}}:=I(1)_{qp^{m_1}}\cdots 
 %           I(s)_{qp^{m_s}},\\
%			 \mathbf{J}(p^b)^{q\underline{n}}&:=J(1)_{p^b}^{[qp^{n_1}]}\cdots J(r)_{p^b}^{[qp^{n_r}]}, 	\mathbf{J}_{p^a}(p^b)^{q\underline{n}}:=J(1)_{p^a,p^b}^{[qp^{n_1}]}\cdots J(r)_{p^a,p^b}^{[qp^{n_r}]}, \mathbf{J}_{p^a,q\underline{n}}:=J(1)_{p^a,qp^{n_1}}\cdots J(r)_{p^a,qp^{n_r}}, \\
%			\mathbf{I}(p^b)^{q\underline{m}}&:=I(1)_{p^b}^{[qp^{m_1}]}\cdots I(s)_{p^b}^{[qp^{m_s}]},\mathbf{I}_{p^a}(p^b)^{q\underline{m}}:=I(1)_{p^a,p^b}^{[qp^{m_1}]}\cdots I(s)_{p^a,p^b}^{[qp^{m_s}]},\mathbf{I}_{p^a,q\underline{m}}:=I(1)_{p^a,qp^{m_1}}\cdots I(s)_{p^a,qp^{m_s}}.
%		\end{split}
%	\end{align*}
{\footnotesize
\begin{align*}
		\begin{split}
			\mathbf{J}_{q\underline{n}}&:=J(1)_{qp^{n_1}}\cdots J(r)_{qp^{n_r}}, \mathbf{I}_{q\underline{m}}:=I(1)_{qp^{m_1}}\cdots 
            I(s)_{qp^{m_s}},\\
			 \mathbf{J}(p^b)^{q\underline{n}}&:=J(1)_{p^b}^{[qp^{n_1}]}\cdots J(r)_{p^b}^{[qp^{n_r}]}, 	\mathbf{J}_{p^a}(p^b)^{q\underline{n}}:=J(1)_{p^a,p^b}^{[qp^{n_1}]}\cdots J(r)_{p^a,p^b}^{[qp^{n_r}]}, \mathbf{J}_{p^a,q\underline{n}}:=J(1)_{p^a,qp^{n_1}}\cdots J(r)_{p^a,qp^{n_r}}, \\
			\mathbf{I}(p^b)^{q\underline{m}}&:=I(1)_{p^b}^{[qp^{m_1}]}\cdots I(s)_{p^b}^{[qp^{m_s}]},\mathbf{I}_{p^a}(p^b)^{q\underline{m}}:=I(1)_{p^a,p^b}^{[qp^{m_1}]}\cdots I(s)_{p^a,p^b}^{[qp^{m_s}]},\mathbf{I}_{p^a,q\underline{m}}:=I(1)_{p^a,qp^{m_1}}\cdots I(s)_{p^a,qp^{m_s}}.
		\end{split}
	\end{align*}
 \par}
\begin{comment}
\begin{tabular}{ll}
  $\mathbf{J}_{q\underline{n}}:=J(1)_{qp^{n_1}}\cdots J(r)_{qp^{n_r}}$,   &  $\mathbf{I}_{q\underline{m}}:=I(1)_{qp^{m_1}}\cdots 
            I(s)_{qp^{m_s}}$, \\
   $\mathbf{J}(p^b)^{q\underline{n}}:=J(1)_{p^b}^{[qp^{n_1}]}\cdots J(r)_{p^b}^{[qp^{n_r}]}$,  & $\mathbf{I}(p^b)^{q\underline{m}}:=I(1)_{p^b}^{[qp^{m_1}]}\cdots I(s)_{p^b}^{[qp^{m_s}]}$,\\
   $\mathbf{J}_{p^a}(p^b)^{q\underline{n}}:=J(1)_{p^a,p^b}^{[qp^{n_1}]}\cdots J(r)_{p^a,p^b}^{[qp^{n_r}]}$, & $\mathbf{I}_{p^a}(p^b)^{q\underline{m}}:=I(1)_{p^a,p^b}^{[qp^{m_1}]}\cdots I(s)_{p^a,p^b}^{[qp^{m_s}]}$, \\
   $\mathbf{J}_{p^a,q\underline{n}}:=J(1)_{p^a,qp^{n_1}}\cdots J(r)_{p^a,qp^{n_r}}$, & $\mathbf{I}_{p^a,q\underline{m}}:=I(1)_{p^a,qp^{m_1}}\cdots I(s)_{p^a,qp^{m_s}}.$
\end{tabular}
\par
\vspace{0.5em}
\end{comment}
We also define the pairs of $p$-families
\begin{equation*}
	\begin{split}
        (\mathscr{J}_{\underline{n}},\mathscr{H}_{\underline{m},\underline{n}})&:=\left(\{\mathbf{J}_{q\underline{n}}\}_{q=1}^{\infty}, \{\mathbf{I}_{q\underline{m}}\mathbf{J}_{q\underline{n}}\}_{q=1}^{\infty}\right), \\
        (\mathscr{J}(p^b)_{\underline{n}},\mathscr{H}(p^b)_{\underline{m},\underline{n}})&:=\left(\{\mathbf{J}(p^b)^{q\underline{n}}\}_{q=1}^{\infty}, \{\mathbf{I}(p^b)^{q\underline{m}}\mathbf{J}(p^b)^{q\underline{n}}\}_{q=1}^{\infty}\right)\text{\space and\space}\\
        (\mathscr{J}_{p^a,\underline{n}},\mathscr{H}_{p^a,\underline{m},\underline{n}})&:=\left(\{\mathbf{J}_{p^a,q\underline{n}}\}_{q=1}^{\infty}, \{\mathbf{I}_{p^a,q\underline{m}}\mathbf{J}_{p^a,q\underline{n}}\}_{q=1}^{\infty}\right).
	\end{split}
\end{equation*}
\end{stp}

%\textcolor{red}{Do we want to include $\mathscr{J}_{p^a}$ here? We use this in the proof of Theorem 3.9.}

In the rest of this section, we further assume that each $(\mathscr{J}_{\underline{n}},\mathscr{H}_{\underline{m},\underline{n}})$ has linear growth, and that if $c_{\underline{m},\underline{n}}:=c(\mathscr{J}_{\underline{n}},\mathscr{H}_{\underline{m},\underline{n}})$, then $c(\mathscr{J}(p^b)_{\underline{n}},\mathscr{H}(p^b)_{\underline{m},\underline{n}})=c_{\underline{m},\underline{n}}.p^b$ for every $b\in \mathbb{N}$. It is worth pointing out that under the assumptions in Setup \ref{pfamiliessetup} that $R$ is a complete local domain, by \cite[Corollary 3.7]{DJ18}, the existence of an OK valuation as in Setup \ref{setuplocaldomainwithOKvaluation} is guaranteed. We also note that the linear growth condition automatically holds when $\mathcal{J}(1), \mathcal{J}(2), \ldots,\mathcal{J}(r)$ are all given by the ring $R$. Furthermore, if, for instance, those $p$-families $\mathcal{I}$ and $\mathcal{J}$ are given by Frobenius powers of ideals, then the condition $c(\mathscr{J}(p^b)_{\underline{n}},\mathscr{H}(p^b)_{\underline{m},\underline{n}})=c_{\underline{m},\underline{n}}.p^b$ is satisfied. \par
\vspace{0.5em}
The following result provides an approximation of the volume of (truncations) of the $p$-body of $p$-families via their shifted families.

\begin{lem}\label{lemmaaboutlimofvolumeJ(p)andH(p)}
	Assume Setup \ref{multipfamiliessetup}. Moreover, assume that $\mathcal{J}(i)$ and $\mathcal{I}(j)$ are all of finite type for $1\leq i\leq r$ and $1\leq j\leq s$. For every $\underline{m}=(m_1,\ldots,m_s)\in \mathbb{N}^s$ and $\underline{n}=(n_1,\ldots n_r)\in \mathbb{N}^r$ we have
	\begin{equation*}
		\lim\limits_{b\to \infty}\lim\limits_{q\to \infty}\dfrac{\#\left(\left[\vartheta^{(h)}\left(\mathscr{J}(p^b)_{\underline{n}}\right)\right]_{q}\cap p^bqH\right)}{p^{bd}q^{d}}
        =
        \mathrm{Vol}_{\mathbb{R}^d}\left(\Delta\left(S, \vartheta^{(h)}\left(\mathscr{J}_{\underline{n}}\right)\right)\cap H\right)\text{\space and \space}
	\end{equation*} 
	\begin{equation*}
		\lim\limits_{b\to \infty}\lim\limits_{q\to \infty}\dfrac{\#\left(\left[\vartheta^{(h)}\left(\mathscr{H}(p^b)_{\underline{m},\underline{n}}\right)\right]_{q}\cap p^bqH\right)}{p^{bd}q^{d}}
        =
        \mathrm{Vol}_{\mathbb{R}^d}\left(\Delta\left(S, \vartheta^{(h)}\left(\mathscr{H}_{\underline{m},\underline{n}}\right)\right)\cap H\right)
	\end{equation*} 
for all $1\leq h\leq [\mathbbm{k}_{\vartheta}:\mathbbm{k}]$ and any truncating half-sapce $H$, where, for instance, $\left[\vartheta^{(h)}\left(\mathscr{A}_{\underline{m},\underline{n}}\right)\right]_{p^b}$ denotes the $p^b$-th element of $\vartheta^{(h)}\left(\mathscr{A}_{\underline{m},\underline{n}}\right)$.
\end{lem}

\begin{proof}
	We fix $\underline{m}=(m_1,\ldots,m_s)\in \mathbb{N}^s$ and $\underline{n}=(n_1,\ldots n_r)\in \mathbb{N}^r$ and consider the following $p$-families:
	\begin{equation*}
		\left(\mathscr{A}_{\underline{m},\underline{n}},\mathscr{B}_{\underline{m},\underline{n}}\right)=\left(\{\mathbf{J}(p^b)^{\underline{n}}\}_{b=1}^{\infty}, \{\mathbf{I}(p^b)^{\underline{m}}\mathbf{J}(p^b)^{\underline{n}}\}_{b=1}^{\infty}\right).
	\end{equation*}
Note that $c(\mathscr{A}_{\underline{m},\underline{n}},\mathscr{B}_{\underline{m},\underline{n}})=c_{\underline{m},\underline{n}}$, then 
\begin{equation*}
	q\left[\vartheta^{(h)}\left(\mathscr{A}_{\underline{m},\underline{n}}\right)\right]_{p^b}+S      \subseteq 
    \left[\vartheta^{(h)}\left(\mathscr{J}(p^b)_{\underline{n}}\right)\right]_{q}+S     \subseteq  \left[\vartheta^{(h)}\left(\mathscr{J}(p^b)_{\underline{n}}\right)\right]_{q}        \subseteq \left[\vartheta^{(h)}\left(\mathscr{A}_{\underline{m},\underline{n}}\right)\right]_{p^bq},
\end{equation*}	
for all $q=p^e,e\in \mathbb{N}$, $b, h\in \mathbb{N}$.
More precisely, the first and second containment hold since $\vartheta^{(h)}\left(\mathscr{A}_{\underline{m},\underline{n}}\right)$ and $\vartheta^{(h)}\left(\mathscr{J}(p^b)_{\underline{n}}\right)$ are $p$-system of ideals in $S$, and the last containment follows by the fact that $\left[\mathscr{J}(p^b)_{\underline{n}}\right]_{q}        \subseteq \left[\mathscr{A}_{\underline{m},\underline{n}}\right]_{p^bq}$ as ideals of $R$. Similarly,
\begin{equation*}	q\left[\vartheta^{(h)}\left(\mathscr{B}_{\underline{m},\underline{n}}\right)\right]_{p^b}+S\subseteq \left[\vartheta^{(h)}\left(\mathscr{H}(p^b)_{\underline{m},\underline{n}}\right)\right]_{q}\subseteq \left[\vartheta^{(h)}\left(\mathscr{B}_{\underline{m},\underline{n}}\right)\right]_{p^bq}.
\end{equation*}
for all $q=p^e,e\in \mathbb{N}$, $b, h\in \mathbb{N}$. Thus, by Theorem \ref{limit=volumeforp-familiesthm} and Theorem \ref{Fujitatypeapproximationtheorem}, for a fixed $\varepsilon\in \mathbb{R}_{>0}$, there exists $b_0\in \mathbb{N}$ such that if $b\geq b_0$ we have
\begin{equation}\label{sosanh3}\begin{split}
		\mathrm{Vol}_{\mathbb{R}^d}(\Delta(S,\vartheta^{(h)}(\mathscr{A}_{\underline{m},\underline{n}}))\cap H)&=\lim\limits_{q\to \infty}\dfrac{\#\left(\left[\vartheta^{(h)}\left(\mathscr{A}_{\underline{m},\underline{n}}\right)\right]_{p^bq}\cap p^bqH\right)}{p^{bd}q^{d}}\\
		&\geq\lim\limits_{q\to \infty}\dfrac{\#\left(\left[\vartheta^{(h)}\left(\mathscr{J}(p^b)_{\underline{n}}\right)\right]_{q}\cap p^bqH\right)}{p^{bd}q^{d}}\\
		&\geq\lim\limits_{q\to \infty}\dfrac{\#\left(\left(q\left[\vartheta^{(h)}\left(\mathscr{A}_{\underline{m},\underline{n}}\right)\right]_{p^b}+S\right)\cap p^bqH\right)}{p^{bd}q^{d}}\\
		&\geq \mathrm{Vol}_{\mathbb{R}^d}(\Delta(S,\vartheta^{(h)}(\mathscr{A}_{\underline{m},\underline{n}}))\cap H)-\varepsilon.
	\end{split}
\end{equation} 
and similarly,
{\footnotesize
\begin{equation}\label{sosanh4}\begin{split}
 		\mathrm{Vol}_{\mathbb{R}^d}(\Delta(S,\vartheta^{(h)}(\mathscr{B}_{\underline{m},\underline{n}}))\cap H) \geq\lim\limits_{q\to \infty}\dfrac{\#\left(\left[\vartheta^{(h)}\left(\mathscr{H}(p^b)_{\underline{m},\underline{n}}\right)\right]_{q}\cap p^bqH\right)}{p^{bd}q^{d}} \geq \mathrm{Vol}_{\mathbb{R}^d}(\Delta(S,\vartheta^{(h)}(\mathscr{B}_{\underline{m},\underline{n}}))\cap H)-\varepsilon,
 	\end{split}
\end{equation}}
for every $h\in \mathbb{Z}_{>0}$ and any truncating half-space $H$. Then, since $\epsilon$ is arbitrary, by taking the limits as $b\rightarrow \infty$, from Inequalities (\ref{sosanh3}) and (\ref{sosanh4}), and Theorem \ref{limit=volumeforp-familiesthm} we deduce that 
\begin{equation}\label{limitJA}
	\lim\limits_{b\to \infty}\lim\limits_{q\to \infty}\dfrac{\#\left(\left[\vartheta^{(h)}\left(\mathscr{J}(p^b)_{\underline{n}}\right)\right]_{q}\cap p^bqH\right)}{p^{bd}q^{d}}=\mathrm{Vol}_{\mathbb{R}^d}(\Delta(S,\vartheta^{(h)}(\mathscr{A}_{\underline{m},\underline{n}}))\cap H) 
\end{equation}
and 
\begin{equation}\label{limitHB}
	\lim\limits_{b\to \infty}\lim\limits_{q\to \infty}\dfrac{\#\left(\left[\vartheta^{(h)}\left(\mathscr{H}(p^b)_{\underline{m},\underline{n}}\right)\right]_{q}\cap p^bqH\right)}{p^{bd}q^{d}}=\mathrm{Vol}_{\mathbb{R}^d}(\Delta(S,\vartheta^{(h)}(\mathscr{B}_{\underline{m},\underline{n}}))\cap H) 
\end{equation}

Now, by the assumption of finite type, there exists a power of $p$, say $g$, such that 
\begin{equation*}
	J(i)_g^{[p^n]}=J(i)_{p^ng} \text{\space and \space} I(j)_g^{[p^n]}=I(j)_{p^ng} \text{\space for all \space} n\geq 0, 1\leq i\leq r, 1\leq j\leq s.
\end{equation*}  
Therefore,
\begin{equation*}
	\mathbf{J}(qg)^{\underline{n}}=\mathbf{J}_{qg\underline{n}} \text{\space and \space} 	\mathbf{I}(qg)^{\underline{m}}\mathbf{J}(qg)^{\underline{n}}=\mathbf{I}_{qg\underline{m}}\mathbf{J}_{qg\underline{n}} 
\end{equation*}
for all $q=p^e, e\in \mathbb{N}$. Hence, by Theorem \ref{limit=volumeforp-familiesthm} we obtain 
\begin{equation}\label{limitsequal1}\begin{split}
		\mathrm{Vol}_{\mathbb{R}^d}(\Delta(S,\vartheta^{(h)}(\mathscr{J}_{\underline{n}}))\cap H)&=\lim\limits_{q\to \infty}\dfrac{\#\left(\left[\vartheta^{(h)}\left(\mathscr{J}_{\underline{n}}\right)\right]_{qg}\cap qgH\right)}{q^{d}g^{d}}\\
		&=\lim\limits_{q\to \infty}\dfrac{\#\left(\left[\vartheta^{(h)}\left(\mathscr{A}_{\underline{m},\underline{n}}\right)\right]_{qg}\cap qgH\right)}{q^{d}g^{d}}\\
		&= \mathrm{Vol}_{\mathbb{R}^d}(\Delta(S,\vartheta^{(h)}(\mathscr{A}_{\underline{m},\underline{n}}))\cap H),
	\end{split}
\end{equation}
and similarly,
\begin{equation}\label{limitsequal2}\begin{split}
		\mathrm{Vol}_{\mathbb{R}^d}(\Delta(S,\vartheta^{(h)}(\mathscr{H}_{\underline{m},\underline{n}}))\cap H) = \mathrm{Vol}_{\mathbb{R}^d}(\Delta(S,\vartheta^{(h)}(\mathscr{B}_{\underline{m},\underline{n}}))\cap H)
	\end{split}
\end{equation} 
The conclusion follows from combining (\ref{limitJA}), (\ref{limitHB}), (\ref{limitsequal1}), and (\ref{limitsequal2}).
\end{proof}

One can also interpret the left-hand side of the formulae in Lemma \ref{lemmaaboutlimofvolumeJ(p)andH(p)} as the volume of certain analog of $p$-body. The following remark illustrates that fact, which is in the spirit of \cite[Lemma 2.8]{YJ22} but also provides a correction to that result.

\begin{rem}
For each $b\in \mathbb{N}$, similar to Definition \ref{def.p-system}, one can define
a collection of subsets $T_{\bullet}=\{T_{qp^b}\}_{q=1}^{\infty}$ of a semigroup $S$ indexed by $q=p^e$ to be a $p^b$-shifted $p$-system if it satisfies the following conditions 
	\begin{enumerate}
		\item $T_{qp^b}$ is an ideal of $S$ for all $q$,
		\item $pT_{qp^b}\subseteq T_{pqp^b}$ for all $q$.
	\end{enumerate}
Similarly, one can define a $p^b$-shifted $p$-family of ideals $\mathcal{I}=\{I_{qp^b}\}_{q=1}^{\infty}$ indexed by powers of $p$ satisfying $I_{qp^b}^{[p]}\subseteq I_{pqp^b}$ for all $q$. Note that $\left\{\mathscr{J}(p^b)_{\underline{n}} \right\}$ is a $p^b$-shifted $p$-family of ideals in $R$, thus, $\vartheta^{(h)}\left(\mathscr{J}(p^b)_{\underline{n}}\right)$ is a $p^b$-shifted $p$-system in $S$. One can define a $p$-body of this system as 
$$\Delta\left(S, \vartheta^{(h)}\left(\mathscr{J}(p^b)_{\underline{n}}\right)\right):=\bigcup_{q=1}^{\infty}\dfrac{1}{qp^b}\left[\vartheta^{(h)}\left(\mathscr{J}(p^b)_{\underline{n}}\right)\right]_q+\mathrm{Cone}(S).$$
By a version of Theorem \ref{limit=volumetheorem}, we have
$$\mathrm{Vol}_{\mathbb{R}^d}\left(\Delta\left(S, \vartheta^{(h)}\left(\mathscr{J}(p^b)_{\underline{n}}\right)\right)\right) = \lim\limits_{q\to \infty}\dfrac{\#\left(\left[\vartheta^{(h)}\left(\mathscr{J}(p^b)_{\underline{n}}\right)\right]_{q}\cap p^bqH\right)}{p^{bd}q^{d}}$$
for any truncating half-space $H$. Thus, Lemma \ref{lemmaaboutlimofvolumeJ(p)andH(p)} can be rewritten in terms of volumes as
\[
   \lim\limits_{b\to \infty} \mathrm{Vol}_{\mathbb{R}^d}\left(\Delta\left(S, \vartheta^{(h)}\left(\mathscr{J}(p^b)_{\underline{n}}\right)\right)\right) = \mathrm{Vol}_{\mathbb{R}^d}\left(\Delta\left(S, \vartheta^{(h)}\left(\mathscr{J}_{\underline{n}}\right)\right)\right)
\]
and 
\[
   \lim\limits_{b\to \infty} \mathrm{Vol}_{\mathbb{R}^d}\left(\Delta\left(S, \vartheta^{(h)}\left(\mathscr{H}(p^b)_{\underline{m},\underline{n}}\right)\right)\right) = \mathrm{Vol}_{\mathbb{R}^d}\left(\Delta\left(S, \vartheta^{(h)}\left(\mathscr{H}_{\underline{m},\underline{n}}\right)\right)\right)
\]

\end{rem}

We are now ready to prove our main result in this section.

\begin{thm}\label{maintool}
	Assume Setup \ref{multipfamiliessetup}. For $\underline{n}=(n_1,\ldots,n_r)\in \mathbb{N}^r,\underline{m}=(m_1,\ldots,m_s)\in \mathbb{N}^s$, and $b\in \mathbb{N}$, we have that the following limits exist and are equal 
	\begin{equation*}
		\lim\limits_{b\to \infty}\lim\limits_{q\to \infty}\dfrac{\lambda(\mathbf{J}(p^b)^{q\underline{n}}/\mathbf{I}(p^b)^{q\underline{m}}\mathbf{J}(p^b)^{q\underline{n}})}{p^{bd}q^d}=\lim\limits_{q\to \infty}\dfrac{\lambda(\mathbf{J}_{q\underline{n}}/\mathbf{I}_{q\underline{m}}\mathbf{J}_{q\underline{n}})}{q^d}.
	\end{equation*}
\end{thm}

\begin{proof}
	Fix $a\in \mathbb{Z}_{>0}$, and note that $I(i)_{p^a,q}\subseteq I(i)_q$ and $J(i)_{p^a,q}\subseteq J(i)_q$ for all $1\leq i\leq s$ and $1\leq j\leq r$, and $q=p^e,e\in \mathbb{N}$. We fix $\underline{n}=(n_1,\ldots,n_r)\in \mathbb{N}^r,\underline{m}=(m_1,\ldots,m_s)\in \mathbb{N}^s$, and for simplicity of notation, we set 
	\begin{equation*}
		(\mathscr{J},\mathscr{H}):=(\mathscr{J}_{\underline{n}},\mathscr{H}_{\underline{m},\underline{n}}),	(\mathscr{J}(p^b),\mathscr{H}(p^b)):=(\mathscr{J}(p^b)_{\underline{n}},\mathscr{H}(p^b)_{\underline{m},\underline{n}}), \text{\space and \space} c=c(\mathscr{J},\mathscr{H}). 
	\end{equation*} 

    As in Lemma \ref{lengthcomputedbyimagesofMunderv(h)}, for the point $\textbf{v}$ in Definition \ref{OKrelativedef}, and set $\textbf{w}=c\textbf{v}$, then $R\cap \mathbb{F}_{\geq q\mathbf{w}}=R\cap \mathbb{F}_{\geq cq\mathbf{v}}\subseteq \mathfrak{m}^{cq}$. Since 
    \begin{align*}
        \mathscr{J}\cap \mathfrak{m}^{cq}=\mathscr{H}\cap \mathfrak{m}^{cq}\text{\space and \space}
        \mathscr{J}(p^b)\cap \mathfrak{m}^{cp^bq}=\mathscr{H}(p^b)\cap \mathfrak{m}^{cp^bq},
    \end{align*}
    then 
    \begin{align*}
        \mathscr{J}\cap \mathbb{F}_{\ge q\textbf{w}}=\mathscr{H}\cap \mathbb{F}_{\ge q\textbf{w}}\text{\space and \space}
        \mathscr{J}(p^b)\cap \mathbb{F}_{\ge qp^b\textbf{w}}=\mathscr{H}(p^b)\cap \mathbb{F}_{\ge qp^b\textbf{w}},
    \end{align*}
    Therefore, we can write 
    \[
    \lambda\left(\dfrac{\mathbf{J}_{q\underline{n}}}{\mathbf{I}_{q\underline{m}}\mathbf{J}_{q\underline{n}}}\right) = \lambda\left(\dfrac{\mathbf{J}_{q\underline{n}}}{\mathbf{J}_{q\underline{n}} \cap \mathbb{F}_{\ge q\textbf{w}}}\right) - \lambda\left(\dfrac{\mathbf{I}_{q\underline{m}}\mathbf{J}_{q\underline{n}}}{\mathbf{I}_{q\underline{m}}\mathbf{J}_{q\underline{n}}\cap \mathbb{F}_{\ge q\textbf{w}}}\right), \text{  and that  }
    \]
    \[
     \lambda\left(\dfrac{\mathbf{J}(p^b)^{q\underline{n}}}{\mathbf{I}(p^b)^{q\underline{m}}\mathbf{J}(p^b)^{q\underline{n}}}\right)
     = \lambda\left(\dfrac{\mathbf{J}(p^b)^{q\underline{n}}}{\mathbf{J}(p^b)^{q\underline{n}} \cap \mathbb{F}_{\ge qp^b\textbf{w}}}\right) - \lambda\left(\dfrac{\mathbf{I}(p^b)^{q\underline{m}}\mathbf{J}(p^b)^{q\underline{n}}}{\mathbf{I}(p^b)^{q\underline{m}}\mathbf{J}(p^b)^{q\underline{n}} \cap \mathbb{F}_{\ge qp^b\textbf{w}}}\right).
     \]  
   \par
   \vspace{0.5em}
	Let $a':=a+\max\{n_1,\ldots,n_r,m_1,\ldots,m_s\}$. It is not hard to check that we get the inclusions $\mathscr{J}_{p^a,q}\subseteq \mathbf{J}_{p^{a'},q\underline{n}}\subseteq\mathbf{J}_{q\underline{n}}=\mathscr{J}_{q}$ for all $q=p^e,e\in \mathbb{N}$. In fact, the first inclusion comes from the fact that $J(i)_{p^tp^{n_i}}^{[q/p^t]}\subseteq J(i)_{p^{a'},qp^{n_i}}$ for all $1\le i \le r$, $0\le t \le a$, and all $q$ power of $p$ as each $\mathcal{J}(i)$ is a $p$-family. The second inclusion follows since $J(i)_{p^a,q}\subseteq J(i)_q$ for all $1\leq j\leq r$, and $q=p^e,e\in \mathbb{N}$. Thus, by Lemma \ref{lengthcomputedbyimagesofMunderv(h)} and Theorem \ref{limit=volumeforp-familiesthm}, with the half-space $H$ given by all $\textbf{u}$ such that $\langle \textbf{u}, \textbf{a} \rangle < \langle \textbf{w}, \textbf{a} \rangle$, where $\textbf{a}$ is the vector defining the valuation in Section \ref{section_about_OKvaluation}, we obtain the following 
	\begin{equation}\label{sosanh5}
		\begin{split}
			\sum\limits_{h=1}^{[\mathbbm{k}_{\vartheta}:\mathbbm{k}]}\mathrm{Vol}_{\mathbb{R}^d}(\Delta(S,\vartheta^{(h)}(\mathscr{J}_{p^a}))\cap H)&=\lim\limits_{q\to \infty}\dfrac{1}{q^d} \cdot \lambda\left(\frac{\mathscr{J}_{p^a,q}}{\mathscr{J}_{p^a,q}\cap \mathbb{F}_{\geq q\mathbf{w}}}\right)\\
			&\leq \lim\limits_{q\to \infty}\dfrac{1}{q^d} \cdot\lambda\left(\frac{\mathbf{J}_{p^{a'},q\underline{n}}}{\mathbf{J}_{p^{a'},q\underline{n}}\cap \mathbb{F}_{\geq q\mathbf{w}}}\right)\\
			&\leq \lim\limits_{q\to \infty}\dfrac{1}{q^d} \cdot\lambda\left(\frac{\mathscr{J}_{q}}{\mathscr{J}_{q}\cap \mathbb{F}_{\geq q\mathbf{w}}}\right)\\
			&=\sum\limits_{h=1}^{[\mathbbm{k}_{\vartheta}:\mathbbm{k}]}\mathrm{Vol}_{\mathbb{R}^d}(\Delta(S,\vartheta^{(h)}(\mathscr{J}))\cap H).
		\end{split}
	\end{equation}

Similarly, we also have $\mathscr{H}_{p^a,q}\subseteq \mathbf{I}_{p^{a'},q\underline{m}}\mathbf{J}_{p^{a'},q\underline{n}}\subseteq\mathbf{I}_{q\underline{m}}\mathbf{J}_{q\underline{n}}=\mathscr{H}_{q}$ for all $q=p^e,e\in \mathbb{N}$, thus,

	{\footnotesize\begin{equation}\label{sosanh6}
	\begin{split}
		\sum\limits_{h=1}^{[\mathbbm{k}_{\vartheta}:\mathbbm{k}]}\mathrm{Vol}_{\mathbb{R}^d}(\Delta(S,\vartheta^{(h)}(\mathscr{H}_{p^a}))\cap H) 
        \leq \lim\limits_{q\to \infty}\dfrac{1}{q^d} 
        \cdot\lambda\left(\frac{\mathbf{I}_{p^{a'},q\underline{m}}\mathbf{J}_{p^{a'},q\underline{n}}}{\mathbf{I}_{p^{a'},q\underline{m}}\mathbf{J}_{p^{a'},q\underline{n}}\cap \mathbb{F}_{\geq q\mathbf{w}}}\right)
        \leq \sum\limits_{h=1}^{[\mathbbm{k}_{\vartheta}:\mathbbm{k}]}\mathrm{Vol}_{\mathbb{R}^d}(\Delta(S,\vartheta^{(h)}(\mathscr{H}))\cap H).
	\end{split}
\end{equation}}
    
    Therefore, using Theorem \ref{thm.volume=volumeoftruncation} (apply the limit as $a\rightarrow \infty$ to the left most term of each inequality (\ref{sosanh5}), (\ref{sosanh6})), we have for any fixed $\varepsilon>0$, there exists $a:=a(\varepsilon)\in \mathbb{Z}_{>0}$ such that 
{\footnotesize\begin{equation}\label{sosanh7}
	\sum\limits_{h=1}^{[\mathbbm{k}_{\vartheta}:\mathbbm{k}]}\mathrm{Vol}_{\mathbb{R}^d}(\Delta(S,\vartheta^{(h)}(\mathscr{J}))\cap H)\geq \lim\limits_{q\to \infty}\dfrac{1}{q^d} \cdot\lambda\left(\frac{\mathbf{J}_{p^{a'},q\underline{n}}}{\mathbf{J}_{p^{a'},q\underline{n}}\cap \mathbb{F}_{\geq q\mathbf{w}}}\right)\geq \sum\limits_{h=1}^{[\mathbbm{k}_{\vartheta}:\mathbbm{k}]}\mathrm{Vol}_{\mathbb{R}^d}(\Delta(S,\vartheta^{(h)}(\mathscr{J}))\cap H)-\dfrac{\varepsilon}{2}, \text{ and}
\end{equation}
\begin{equation}\label{sosanh8}	
\sum\limits_{h=1}^{[\mathbbm{k}_{\vartheta}:\mathbbm{k}]}\mathrm{Vol}_{\mathbb{R}^d}(\Delta(S,\vartheta^{(h)}(\mathscr{H}))\cap H)\geq \lim\limits_{q\to \infty}\dfrac{1}{q^d} \cdot\lambda\left(\frac{\mathbf{I}_{p^{a'},q\underline{m}}\mathbf{J}_{p^{a'},q\underline{n}}}{\mathbf{I}_{p^{a'},q\underline{m}}\mathbf{J}_{p^{a'},q\underline{n}}\cap \mathbb{F}_{\geq q\mathbf{w}}}\right)\geq \sum\limits_{h=1}^{[\mathbbm{k}_{\vartheta}:\mathbbm{k}]}\mathrm{Vol}_{\mathbb{R}^d}(\Delta(S,\vartheta^{(h)}(\mathscr{H}))\cap H)-\dfrac{\varepsilon}{2}.
\end{equation}}

On the other hand, as each $\mathcal{I}(i)_{p^{a'}}$ and each $\mathcal{J}(j)_{p^{a'}}$ is a $p$-family of finite type, thus by Lemma \ref{lemmaaboutlimofvolumeJ(p)andH(p)}, we have 
\begin{equation}\label{estimateJ_(p^(a'))p^bwithJ_(p^(a'))}
    \begin{split}
        \lim\limits_{b\to \infty}\lim\limits_{q\to \infty} \dfrac{1}{p^{bd}q^d}\cdot \lambda\left(\frac{\mathbf{J}_{p^{a'}}(p^b)^{q\underline{n}}}{\mathbf{J}_{p^{a'}}(p^b)^{q\underline{n}}\cap \mathbb{F}_{\geq qp^b\mathbf{w}}}\right)
        &=\sum\limits_{h=1}^{[\mathbbm{k}_{\vartheta}:\mathbbm{k}]}\lim\limits_{b\to \infty}\lim\limits_{q\to \infty}\dfrac{\#\left(\vartheta^{(h)}\left(\mathbf{J}_{p^{a'}}(p^b)^{q\underline{n}}\right)\cap p^bqH\right)}{p^{bd}q^{d}}\\
        &= \sum\limits_{h=1}^{[\mathbbm{k}_{\vartheta}:\mathbbm{k}]} \mathrm{Vol}_{\mathbb{R}^d}\left(\Delta\left(S, \vartheta^{(h)}\left(\mathscr{J}_{p^{a'},\underline{n}}\right)\right)\cap H\right)\\
        &=\lim\limits_{q\to \infty}\sum\limits_{h=1}^{[\mathbbm{k}_{\vartheta}:\mathbbm{k}]}\dfrac{\#\left(\vartheta^{(h)}\left(\mathbf{J}_{p^{a'},q\underline{n}}\right)\cap qH\right)}{q^{d}} \\
        &=\lim\limits_{q\to \infty} \dfrac{1}{q^d}\cdot \lambda\left(\frac{\mathbf{J}_{p^{a'},q\underline{n}}}{\mathbf{J}_{p^{a'},q\underline{n}}\cap \mathbb{F}_{\geq q\mathbf{w}}}\right),
    \end{split}  
\end{equation}
and, similarly,
\begin{equation}\label{estimateH_(p^(a'))p^bwithH_(p^(a'))}
     \lim\limits_{b\to \infty}\lim\limits_{q\to \infty} \dfrac{1}{p^{bd}q^d}\cdot \lambda\left(\frac{\mathbf{I}_{p^{a'}}(p^b)^{q\underline{m}}\mathbf{J}_{p^{a'}}(p^b)^{q\underline{n}}}{\mathbf{I}_{p^{a'}}(p^b)^{q\underline{m}}\mathbf{J}_{p^{a'}}(p^b)^{q\underline{n}}\cap \mathbb{F}_{\geq qp^b\mathbf{w}}}\right)=\lim\limits_{q\to \infty} \dfrac{1}{q^d}\cdot \lambda\left(\frac{\mathbf{I}_{p^{a'},q\underline{m}}\mathbf{J}_{p^{a'},q\underline{n}}}{\mathbf{I}_{p^{a'},q\underline{m}}\mathbf{J}_{p^{a'},q\underline{n}}\cap \mathbb{F}_{\geq q\mathbf{w}}}\right).
\end{equation}
On the other hand, for every $b\in \mathbb{Z}_{>0}$, we have the inclusions 
\begin{equation*}
	\mathbf{J}_{p^{a'}}(p^b)^{q\underline{n}}\subseteq \mathscr{J}(p^b)_{q}\subseteq \mathscr{J}_{p^bq} \text{\space and \space} \mathbf{I}_{p^{a'}}(p^b)^{q\underline{m}}\mathbf{J}_{p^{a'}}(p^b)^{q\underline{n}}\subseteq \mathscr{H}(p^b)_{q}\subseteq \mathscr{H}_{p^bq}. 
\end{equation*}
Therefore by Equalities (\ref{estimateJ_(p^(a'))p^bwithJ_(p^(a'))}), (\ref{estimateH_(p^(a'))p^bwithH_(p^(a'))}), and Lemma \ref{lengthcomputedbyimagesofMunderv(h)} there exists $b_0:=b_0(a)$ such that if $b\geq b_0$ we have 
\begin{equation}\label{sosanh9}
	\begin{split}
\sum\limits_{h=1}^{[\mathbbm{k}_{\vartheta}:\mathbbm{k}]}\mathrm{Vol}_{\mathbb{R}^d}(\Delta(S,\vartheta^{(h)}(\mathscr{J}))\cap H)&=\lim\limits_{q\to \infty}\sum\limits_{h=1}^{[\mathbbm{k}_{\vartheta}:\mathbbm{k}]}\dfrac{\#\left(\left[\vartheta^{(h)}(\mathscr{I})\right]_{p^{b}q}\cap p^bqH \right)}{p^{bd}q^{d}}\\
    &\geq \lim\limits_{q\to \infty}\sum\limits_{h=1}^{[\mathbbm{k}_{\vartheta}:\mathbbm{k}]}\dfrac{\#\left(\left[\vartheta^{(h)}(\mathscr{I}(p^b))\right]_{q}\cap p^bqH \right)}{p^{bd}q^{d}}\\
    &\geq \lim\limits_{q\to \infty}\sum\limits_{h=1}^{[\mathbbm{k}_{\vartheta}:\mathbbm{k}]}\dfrac{\#\left(\vartheta^{(h)}(\mathbf{J}_{p^{a'}}(p^b)^{q\underline{n}})\cap p^bqH \right)}{p^{bd}q^{d}}\\
	&= \lim\limits_{q\to \infty}\dfrac{1}{p^{bd}q^d}\cdot \lambda\left(\frac{\mathbf{J}_{p^{a'}}(p^b)^{q\underline{n}}}{\mathbf{J}_{p^{a'}}(p^b)^{q\underline{n}}\cap \mathbb{F}_{\geq qp^b\mathbf{w}}}\right)\\
	&\geq \lim\limits_{q\to \infty}\dfrac{1}{q^d}\cdot \lambda\left(\frac{\mathbf{J}_{p^{a'},q\underline{n}}}{\mathbf{J}_{p^{a'},q\underline{n}}\cap \mathbb{F}_{\geq q\mathbf{w}}}\right)-\dfrac{\varepsilon}{2},
	\end{split}
\end{equation}
and similarly
\begin{equation}\label{sosanh10}
	\begin{split}
	\sum\limits_{h=1}^{[\mathbbm{k}_{\vartheta}:\mathbbm{k}]}\mathrm{Vol}_{\mathbb{R}^d}(\Delta(S,\vartheta^{(h)}(\mathscr{H}))\cap H)&\geq \lim\limits_{q\to \infty}\dfrac{1}{p^{bd}q^d}\cdot \lambda\left(\frac{\mathbf{I}_{p^{a'}}(p^b)^{q\underline{m}}\mathbf{J}_{p^{a'}}(p^b)^{q\underline{n}}}{\mathbf{I}_{p^{a'}}(p^b)^{q\underline{m}}\mathbf{J}_{p^{a'}}(p^b)^{q\underline{n}}\cap \mathbb{F}_{\geq qp^b\mathbf{w}}}\right)\\
	&\geq \lim\limits_{q\to \infty}\dfrac{1}{q^d}\cdot \lambda\left(\frac{\mathbf{I}_{p^{a'},q\underline{m}}\mathbf{J}_{p^{a'},q\underline{n}}}{\mathbf{I}_{p^{a'},q\underline{m}}\mathbf{J}_{p^{a'},q\underline{n}}\cap \mathbb{F}_{\geq q\mathbf{w}}}\right)-\dfrac{\varepsilon}{2}.
	\end{split}
\end{equation} 
The result now follows by combining the Inequalities (\ref{sosanh7}) with (\ref{sosanh9}), and (\ref{sosanh8}) with (\ref{sosanh10}) as $\epsilon$ is arbitrary.
\end{proof}

\begin{thm}\label{thm.maintool}
   Let $(R,\mathfrak{m},\mathbbm{k})$ be a Noetherian local ring of dimension $d$ and characteristic $p>0$ with perfect residue field $\mathbbm{k}$ such that $\dim (\mathrm{N}(\hat{R}))<d$; here $\mathrm{N}(\hat{R})$ denotes the nilradical of the $\mathfrak{m}$-adic completion $\hat{R}$.  Let $\mathcal{J}(1)=\{J(1)_q\}_{q=1}^{\infty},\ldots,\mathcal{J}(r)=\{J(r)_q\}_{q=1}^{\infty}$ be $p$-families of non-zero ideals, and let $\mathcal{I}(1)=\{I(1)_q\}_{q=1}^{\infty},\ldots,\mathcal{I}(s)=\{I(s)_q\}_{q=1}^{\infty}$ be $\mathfrak{m}$-primary $p$-families of ideals. For $\underline{n}=(n_1,\ldots,n_r)\in \mathbb{N}^r,\underline{m}=(m_1,\ldots,m_s)\in \mathbb{N}^s$, and $b,a\in \mathbb{N}$, we follow the notations in Setup \ref{multipfamiliessetup}. Then, we have that the following limits exist and are equal 
	\begin{equation*}
		\lim\limits_{b\to \infty}\lim\limits_{q\to \infty}\dfrac{\lambda(\mathbf{J}(p^b)^{q\underline{n}}/\mathbf{I}(p^b)^{q\underline{m}}\mathbf{J}(p^b)^{q\underline{n}})}{p^{bd}q^d}=\lim\limits_{q\to \infty}\dfrac{\lambda(\mathbf{J}_{q\underline{n}}/\mathbf{I}_{q\underline{m}}\mathbf{J}_{q\underline{n}})}{q^d}.
	\end{equation*}
\end{thm}

\begin{proof}
    By passing through the $\mathfrak{m}$-adic completion $\hat{R}$, we can assume $R$ is a complete local ring. By the assumption and \cite[Corollary 5.18]{DJ18}, we can also assume $R$ is complete and reduced. Suppose that the minimal primes of $R$ are $\{P_1,\ldots,P_s\}$. Let $R_i=R/P_i$ for every $1\leq i\leq s$. Then each $R_i$ is a complete local domain for $i=1,\ldots,s$. By Lemma \ref{lem.length.over.R.to.sum.length.over.Ri}, we have 
    \begin{equation*}
		\lim\limits_{b\to \infty}\lim\limits_{q\to \infty}\dfrac{\lambda(\mathbf{J}(p^b)^{q\underline{n}}/\mathbf{I}(p^b)^{q\underline{m}}\mathbf{J}(p^b)^{q\underline{n}})}{p^{bd}q^d}=\sum\limits_{i=1}^s\lim\limits_{b\to \infty}\lim\limits_{q\to \infty}\dfrac{\lambda_{R_i}(\mathbf{J}(p^b)^{q\underline{n}}R_i/\mathbf{I}(p^b)^{q\underline{m}}\mathbf{J}(p^b)^{q\underline{n}}R_i)}{p^{bd}q^d},
    \end{equation*}
    and
    \begin{equation*}
	\lim\limits_{q\to \infty}\dfrac{\lambda(\mathbf{J}_{q\underline{n}}/\mathbf{I}_{q\underline{m}}\mathbf{J}_{q\underline{n}})}{q^d}=\sum\limits_{i=1}^{s}\lim\limits_{q\to \infty}\dfrac{\lambda_{R_i}(\mathbf{J}_{q\underline{n}}R_i/\mathbf{I}_{q\underline{m}}\mathbf{J}_{q\underline{n}}R_i)}{q^d}.
	\end{equation*}
 Therefore, combining the two equations above with Theorem \ref{maintool}, we have 
\begin{align*}
    \lim\limits_{b\to \infty}\lim\limits_{q\to \infty}\dfrac{\lambda(\mathbf{J}(p^b)^{q\underline{n}}/\mathbf{I}(p^b)^{q\underline{m}}\mathbf{J}(p^b)^{q\underline{n}})}{p^{bd}q^d}&=\sum\limits_{i=1}^s\lim\limits_{b\to \infty}\lim\limits_{q\to \infty}\dfrac{\lambda_{R_i}(\mathbf{J}(p^b)^{q\underline{n}}R_i/\mathbf{I}(p^b)^{q\underline{m}}\mathbf{J}(p^b)^{q\underline{n}}R_i)}{p^{bd}q^d}\\
    &=\sum\limits_{i=1}^{s}\lim\limits_{q\to \infty}\dfrac{\lambda_{R_i}(\mathbf{J}_{q\underline{n}}R_i/\mathbf{I}_{q\underline{m}}\mathbf{J}_{q\underline{n}}R_i)}{q^d}.\\
    &=\lim\limits_{q\to \infty}\dfrac{\lambda(\mathbf{J}_{q\underline{n}}/\mathbf{I}_{q\underline{m}}\mathbf{J}_{q\underline{n}})}{q^d}.
\end{align*}
\end{proof}

%%%%%%%%%%%%%%%%
%%%%%%%%%%%%%%%%%%%%%%%%%%%%%%%%%%%%%%
\section{Applications}
\label{sec.app}
In this section, we use our technical results from the previous section and results in \cite{Verma91} to obtain a general version of results in \cite{WY01} concerning the length function of multi-$p$-families, their multiplicities, mixed multiplicities, and Hilbert-Kunz multiplicities in a $2$-dimensional Cohen-Macaulay ring. We also show an instance where one can define by showing the existence of a mixed multiplicity version of multi-$p-$families of ideals.

%%%%%%%%%%%%%%
%%%%%%%%%%%%%%%%%%%%%%%%%%%%%%%%%%%%%

\subsection{$p$-families in Cohen-Macaulay local ring of dimension $2$}

\begin{stp}
\label{setupp-familiesfromgraded}
Let $(R,\mathfrak{m},\mathbbm{k})$ be a Noetherian local ring of dimension $d$ and characteristic $p>0$ with perfect residue field $\mathbbm{k}$ such that $\dim (\mathrm{N}(\hat{R}))<d$, where $\mathrm{N}(\hat{R})$ denotes the nilradical of the $\mathfrak{m}$-adic completion $\hat{R}$. Let $\mathcal{I}(1)=\{I(1)_q\}_{q=1}^{\infty},\ldots,\mathcal{I}(s)=\{I(s)_q\}_{q=1}^{\infty}$ be $\mathfrak{m}$-primary $p$-families of ideals such that 
 \begin{equation}\label{containment.condition.with.ordinary.powers}
  I(i)_q^{[p]}\subseteq I(i)_q^p\subseteq I(i)_{qp}\text{\space for every \space} i=1,\ldots,s   
 \end{equation}
 For every $b\in \mathbb{N}$ and $\underline{m}=(m_1,\ldots,m_s)\in \mathbb{N}^s$ we follow the same abbreviations from Setup \ref{multipfamiliessetup}. The sequence $(\mathcal{I}(1),\ldots,\mathcal{I}(s))$ of $p$-families is simply denoted by $\mathcal{I}$.
\end{stp}

 In the setup above, the first containment of (\ref{containment.condition.with.ordinary.powers}) is obvious but the second containment is not true for a general $p$-family. One of the simplest examples is to consider the $p$-families $\mathcal{I}(i)=\{I(i)_q\}_{q=1}^{\infty}$ obtained from picking the elements indexed by powers of $p$ of $\mathfrak{m}$-primary graded families of ideals $\mathfrak{I}(i)=\{I(i)_n\}_{n=1}^{\infty}$ for $i=1,\ldots,s$. The following example shows that there is no obvious way to fill a $p$-family with property (\ref{containment.condition.with.ordinary.powers}) to a graded family.

\begin{exmp}
    Consider $R=\mathbb{F}_p[x,y,z]$. Let $I=(x,y)$ and $J=(z)$ be ideals in $R$. We construct a $p$-family of ideals $\mathcal{K}=\{K_q\}_{q=1}^{\infty}$ in $R$ as follows
    \begin{align*}
        K_{p^{n}}&=K_{p^{n-1}}^p+J^{[p]} \text{\space for all \space} n\in \mathbb{Z}.
    \end{align*}
    It is obvious that $K_q^{[p]}\subseteq K_q^p\subseteq K_{qp}\text{\space for every \space} i=1,\ldots,s$. 
\end{exmp}

Property (\ref{containment.condition.with.ordinary.powers}) also leads us to an interesting question.

\begin{quest}
        For what $p-$family can we add more ideals to make it a graded family? Certainly, we cannot fill in a Frobenius powers family to obtain a graded family.
\end{quest}

\begin{prop}\label{prop.lim.of.multi.equals.to.lim.of.HKmulti}
    Let $\mathcal{I}=\{I_q\}_{q=1}^{\infty}$ be a $p$-family of ideals in $R$ such that $I_q^{[p]}\subseteq I_q^p\subseteq I_{qp}$ for all $q=p^e,e\in \mathbb{N}$. Then
    \begin{equation*}
        \lim\limits_{q\to \infty}\dfrac{e(I_q)}{q^d}=d!\lim\limits_{q\to \infty}\dfrac{e_{HK}(I_q)}{q^d}.
    \end{equation*}
\end{prop}

\begin{proof}
    Since $I_q^{[p^b]}\subseteq I_q^{p^b}\subseteq I_{qp^b}$ for all $b\in \mathbb{N}, q=p^e, e\in \mathbb{N}$, then 
    \begin{equation}\label{compare_length}
        \lambda\left(R/I_{qp^b}\right)\leq \lambda\left(R/I_q^{p^b}\right)\leq \lambda\left(R/I_q^{[p^b]}\right).
    \end{equation}
It follows that 
\[
\liminf\limits_{q\to \infty}\lim\limits_{b\to \infty} \dfrac{\lambda(R/I_{qp^b})}{p^{bd}q^d}\leq\liminf\limits_{q\to \infty}\lim\limits_{b\to \infty} \dfrac{\lambda(R/I_q^{p^b})}{p^{bd}q^d} \text{  and  } \limsup\limits_{q\to \infty}\lim\limits_{b\to \infty} \dfrac{\lambda(R/I_q^{p^b})}{p^{bd}q^d}\leq \limsup\limits_{q\to \infty}\lim\limits_{b\to \infty}\dfrac{\lambda(R/I_q^{[p^b]})}{p^{bd}q^d}.
\]

\begin{comment}
It follows from the first inequality of (\ref{compare_length}) that 
    \begin{equation*}
       \lim\limits_{q\to \infty}\lim\limits_{b\to \infty} \dfrac{\lambda(R/I_{qp^b})}{p^{bd}q^d}=\liminf\limits_{q\to \infty}\lim\limits_{b\to \infty} \dfrac{\lambda(R/I_{qp^b})}{p^{bd}q^d}\leq\liminf\limits_{q\to \infty}\lim\limits_{b\to \infty} \dfrac{\lambda(R/I_q^{p^b})}{p^{bd}q^d}.
    \end{equation*}
    From the second inequality of (\ref{compare_length}), we have 
    \begin{equation*}
       \limsup\limits_{q\to \infty}\lim\limits_{b\to \infty} \dfrac{\lambda(R/I_q^{p^b})}{p^{bd}q^d}\leq \limsup\limits_{q\to \infty}\lim\limits_{b\to \infty}\dfrac{\lambda(R/I_q^{[p^b]})}{p^{bd}q^d}=\lim\limits_{q\to \infty}\lim\limits_{b\to \infty}\dfrac{\lambda(R/I_q^{[p^b]})}{p^{bd}q^d}.
    \end{equation*}
\end{comment}

On the other hand, note that $\lim\limits_{q\to \infty}\lim\limits_{b\to \infty}\dfrac{\lambda(R/I_q^{[p^b]})}{p^{bd}q^d}=\lim\limits_{q\to \infty}\dfrac{e_{HK}(I_q)}{q^d}$. Moreover, by \cite[Theorem 3.2]{SD22}, we have
    \begin{equation*}
        \lim\limits_{q\to \infty}\lim\limits_{b\to \infty} \dfrac{\lambda(R/I_{qp^b})}{p^{bd}q^d}= \lim\limits_{q\to \infty}\dfrac{\lambda(R/I_{q})}{q^d}=\lim\limits_{q\to \infty}\dfrac{e_{HK}(I_q)}{q^d}.
    \end{equation*}

    Hence, we get
    \begin{equation*}
        \lim\limits_{q\to \infty}\dfrac{e_{HK}(I_q)}{q^d}\leq \liminf\limits_{q\to \infty}\lim\limits_{b\to \infty} \dfrac{\lambda(R/I_q^{p^b})}{p^{bd}q^d}\leq\limsup\limits_{q\to \infty}\lim\limits_{b\to \infty} \dfrac{\lambda(R/I_q^{p^b})}{p^{bd}q^d}\leq \lim\limits_{q\to \infty}\dfrac{e_{HK}(I_q)}{q^d}.
    \end{equation*}
    Therefore, the limits
    \begin{equation*}
         \lim\limits_{q\to \infty}\dfrac{e(I_q)}{q^d}= d!\lim\limits_{q\to \infty}\lim\limits_{n\to \infty}\dfrac{\lambda(R/I_q^n)}{n^dq^d}=d!\lim\limits_{q\to \infty}\lim\limits_{b\to \infty}\dfrac{\lambda(R/I_q^{p^b})}{p^{bd}q^d}
    \end{equation*}
    exist and are equal to $d!\lim\limits_{q\to \infty}\dfrac{e_{HK}(I_q)}{q^d}$.
\end{proof}

\begin{quest}
    In general for any $p$-family of ideals $\mathcal{I}=\{I_q\}_{q=1}^{\infty}$, does the limit $\lim\limits_{q\to \infty}\dfrac{e(I_q)}{q^d}$ exist?
\end{quest}

With Theorem \ref{thm.maintool} and the property that $I(i)_q^{[p]}\subseteq I(i)_q^p\subseteq I(i)_{qp}$ for every $i=1,\ldots,s$ and $q=p^e,e\geq 0$, we have the following consequence: 

\begin{thm}\label{limittheoremingradedcase}
	Assume Setup \ref{setupp-familiesfromgraded}. For $\underline{m}=(m_1,\ldots,m_s)\in \mathbb{N}^s$, and $b\in \mathbb{N}$, we have that the following limits exist and are equal 
	\begin{equation*}
		\lim\limits_{b\to \infty}\lim\limits_{q\to \infty}\dfrac{\lambda(R/\mathbf{I}(p^b)^{q\underline{m}})}{p^{bd}q^d}=\lim\limits_{b\to \infty}\lim\limits_{q\to \infty}\dfrac{\lambda(R/I(1)_{p^b}^{[q].p^{m_1}}\cdots I(s)_{p^b}^{[q].p^{m_s}})}{p^{bd}q^d}=\lim\limits_{q\to \infty}\dfrac{\lambda(R/\mathbf{I}_{q\underline{m}})}{q^d},
	\end{equation*}
where $I^{[q].p^{m}}=(I^{[q]})^{p^{m}}=(I^{p^m})^{[q]}$ for any ideal $I\subset R$.
\end{thm}

%\textcolor{red}{Does $I^{[q].p^{m}}$ mean $(I^{[q]})^{p^{m}}$? Is it equal to $(I^{p^m})^{[q]}$?}

Next, let us recall some results from \cite{Verma91}.

\begin{defn}\cite[Definition 2.7]{Verma91}    
Let $(R,\mathfrak{m})$ be a $2$-dimensional CM local ring. Let $I$ and $J$ be $\mathfrak{m}$-primary ideals. If
    \begin{equation*}
     e(I|J) = \lambda(R/IJ) - \lambda(R/I) - \lambda(R/J)
    \end{equation*}
    we write $r(I|J)=0$. Here the mixed multiplicity of $I$ and $J$ can be defined by
    \begin{equation*}
        e(I|J)=\dfrac{1}{2}(e(IJ)-e(I)-e(J)).
    \end{equation*}
\end{defn}

\begin{thm}\label{Vermatheorem}\cite[Proposition 2.9]{Verma91}    
Let $(R,\mathfrak{m})$ be a $2$-dimensional CM local ring. Let $I_1,I_2,\ldots,I_g$ be $\mathfrak{m}$-primary ideals. If $r(I_i|I_j)=0$ for $1\leq i\leq j\leq g$, then for all non-negative integers $r_1,\ldots,r_s$
    \begin{equation*}
        \lambda(R/I_1^{r_1}\ldots I_s^{r_s})=\sum\limits_{i=1}^g\left(e(I_i)\binom{r_i}{2}+\lambda(R/I_i)r_i\right)+\sum\limits_{1\leq i\leq j\leq g}e(I_i|J_i)r_ir_j.
    \end{equation*}
\end{thm}

By using the theorem above, we can describe explicitly the limit $\lim\limits_{q\to \infty}\frac{\lambda(R/\mathbf{I}_{q\underline{m}})}{q^d}$ when $R$ is a $2$-dimensional Cohen-Macaulay local ring with $\dim (\mathrm{N}(\hat{R}))<2$ and perfect residue field.

\begin{thm} \label{maintheoreminCMlocalringdim2}
    Let $(R,\mathfrak{m})$ be a $2$-dimensional CM local ring. Assume Setup \ref{setupp-familiesfromgraded}. For every $b\in \mathbb{N}$, we assume that $r(I(i)_{p^b}|I(j)_{p^b})=0$ for $1\leq i\leq j\leq s$. Then for $\underline{m}=(m_1,\ldots,m_s)\in \mathbb{N}^s$, we have 
    {\footnotesize
    \begin{equation*}
        \lim\limits_{q\to \infty}\dfrac{\lambda(R/\mathbf{I}_{q\underline{m}})}{q^2}=\sum\limits_{i=1}^s\left(\lim\limits_{b\to \infty}\dfrac{e(I(i)_{p^b})}{p^{2b}}\binom{p^{m_i}}{2}+\lim\limits_{b\to \infty}\dfrac{e_{HK}(I(i)_{p^b})}{p^{2b}}p^{m_i}\right)+\sum\limits_{1\leq i\leq j\leq s}\lim\limits_{b\to \infty}\dfrac{e(I(i)_{p^b}|I(j)_{p^b})}{p^{2b}}p^{m_i}p^{m_j}.
    \end{equation*}}
\end{thm}

\begin{proof}
    By Theorem \ref{Vermatheorem}, for each $b\in \mathbb{N}$, we have
    \begin{align*}
        \lambda(R/I(1)_{p^b}^{[q].p^{m_1}}\cdots I(s)_{p^b}^{[q].p^{m_s}})=&\sum\limits_{i=1}^s\left(e(I(i)_{p^b}^{[q]})\binom{p^{m_i}}{2}+\lambda(R/I(i)_{p^b}^{[q]})p^{m_i}\right)\\
        &+\sum\limits_{1\leq i\leq j\leq s}e(I(i)_{p^b}^{[q]}|I(j)_{p^b}^{[q]})p^{m_i}p^{m_j}.
    \end{align*}
    Note that by \cite[Proposition 3.4]{HTW20}, since $(IJ)^{[p]}=I^{[p]}J^{[p]}$, we have that for all $q=p^e$, $$I(1)_{p^b}^{[q].p^{m_1}}\cdots I(s)_{p^b}^{[q].p^{m_s}}=(I(1)_{p^b}^{p^{m_1}}\cdots I(s)_{p^b}^{p^{m_s}})^{[q]}.$$
    Moreover, since $e(I(i)_{p^b}^{[q]})=q^2e(I(i)_{p^b})$ for all $1\leq i\leq s$ (\cite[Theorem 3.2]{D03}), we have
    \begin{equation*}
        \begin{split}
            e(I(i)_{p^b}^{[q]}|I(j)_{p^b}^{[q]})&=\dfrac{1}{2}\left(e(I(i)_{p^b}^{[q]}.I(j)_{p^b}^{[q]})-e(I(i)_{p^b}^{[q]})-e(I(j)_{p^b}^{[q]})\right)\\
            &= \dfrac{q^2}{2}\left(e(I(i)_{p^b}I(j)_{p^b})-e(I(i)_{p^b})-e(I(j)_{p^b})\right)\\
            &=q^2e(I(i)_{p^b}|I(j)_{p^b}).
        \end{split}
    \end{equation*}
    Therefore,
    \begin{align*}
       & \lim\limits_{q\to \infty}\dfrac{\lambda(R/I(1)_{p^b}^{[q].p^{m_1}}\cdots I(s)_{p^b}^{[q].p^{m_s}})}{q^2}\\
       =&\sum\limits_{i=1}^s\left(e(I(i)_{p^b})\binom{p^{m_i}}{2}+e_{HK}(I(i)_{p^b})p^{m_i}\right)+\sum\limits_{1\leq i\leq j\leq s}e(I(i)_{p^b}|I(j)_{p^b})p^{m_i}p^{m_j}.
    \end{align*}
    By Theorem \ref{limittheoremingradedcase} taking the limits give
    \begin{align*}
        &\lim\limits_{q\to \infty}\dfrac{\lambda(R/\mathbf{I}_{q\underline{m}})}{q^2} = \lim\limits_{b\to \infty}\lim\limits_{q\to \infty}\dfrac{\lambda(R/I(1)_{p^b}^{[q].p^{m_1}}\cdots I(s)_{p^b}^{[q].p^{m_s}})}{q^2p^{2b}}\\
       =&\lim\limits_{b\to \infty}\dfrac{1}{p^{2b}}\left[\sum\limits_{i=1}^s\left(e(I(i)_{p^b})\binom{p^{m_i}}{2}+e_{HK}(I(i)_{p^b})p^{m_i}\right)+\sum\limits_{1\leq i\leq j\leq s}e(I(i)_{p^b}|I(j)_{p^b})p^{m_i}p^{m_j}\right]\\
       =&\sum\limits_{i=1}^s\left(\lim\limits_{b\to \infty}\dfrac{e(I(i)_{p^b})}{p^{2b}}\binom{p^{m_i}}{2}+\lim\limits_{b\to \infty}\dfrac{e_{HK}(I(i)_{p^b})}{p^{2b}}p^{m_i}\right)+\sum\limits_{1\leq i\leq j\leq s}\lim\limits_{b\to \infty}\dfrac{e(I(i)_{p^b}|I(j)_{p^b})}{p^{2b}}p^{m_i}p^{m_j}.
       %=&\sum\limits_{i=1}^s\left(e(\mathfrak{I}(i))\binom{p^{m_i}}{2}+\mathrm{Vol}(\mathcal{I}(i))p^{m_i}\right)+\sum\limits_{1\leq i\leq j\leq s}e(\mathfrak{I}(i)|\mathfrak{I}(j))p^{m_i}p^{m_j}.%
    \end{align*}
    Finally, we show that all coefficients of the polynomial above exist. Since $r(I(i)_{p^b}|I(j)_{p^b})=0$ for $1\leq i\leq j\leq s$  we have 
    \begin{equation*}
        e(I(i)_{p^b}|I(j)_{p^b})= \lambda(R/I(i)_{p^b}I(j)_{p^b}) - \lambda(R/I(i)_{p^b}) - \lambda(R/I(j)_{p^b}).
    \end{equation*}
    It follows that for every $1\leq i\leq j\leq s$, the limit 
    \begin{align*}
        \lim\limits_{b\to \infty}\dfrac{e(I(i)_{p^b}|I(j)_{p^b})}{p^{2b}}&=\lim\limits_{b\to \infty}\dfrac{\lambda(R/I(i)_{p^b}I(j)_{p^b})}{p^{2b}}-\lim\limits_{b\to \infty}\dfrac{\lambda(R/I(i)_{p^b})}{p^{2b}}-\lim\limits_{b\to \infty}\dfrac{\lambda(R/I(j)_{p^b})}{p^{2b}}\\
        %&=\mathrm{Vol}(\mathcal{I}(i)\mathcal{I}(j))-\mathrm{Vol}(\mathcal{I}(i))-\mathrm{Vol}(\mathcal{I}(j))
    \end{align*}
    exists by \cite[Theorem 5.15]{DJ18}. Also, by Proposition \ref{prop.lim.of.multi.equals.to.lim.of.HKmulti}, $\lim\limits_{b\to \infty}\dfrac{e(I(i)_{p^b})}{p^{2b}}$ exists.
\end{proof}

\begin{exmp}
\begin{enumerate}
    \item One of the examples that satisfies the theorem above is a pseudo-rational local ring of dimension $2$ with characteristic $p>0$ such that $\dim (\mathrm{N}(\hat{R}))<2$, and all ideals in Setup \ref{setupp-familiesfromgraded} are $\mathfrak{m}$-primary complete ideals, see \cite{Rees81, WY01} for more details. It can be considered as a generalization of a result of Wantanabe and Yoshida \cite[Proposition 4.6]{WY01}. 
    \item In Setup \ref{setupp-familiesfromgraded}, if the $\mathfrak{m}$-primary $p$-families $\mathcal{I}(1)=\{I(1)_q\}_{q=1}^{\infty},\ldots,\mathcal{I}(s)=\{I(s)_q\}_{q=1}^{\infty}$ are subfamilies of $\mathfrak{m}$-primary graded families of ideals $\mathfrak{I}(1)=\{I(1)_n\}_{n=1}^{\infty},\ldots,\mathfrak{I}(s)=\{I(s)_n\}_{n=1}^{\infty}$ respectively, then the formula in Theorem \ref{maintheoreminCMlocalringdim2} coincides with the result about the mixed multiplicities of the given graded families as in \cite[Theorem 3.3 and Corollary 3.5]{YJ22}. 
    %\textcolor{red}{nothing to do with Verma or Wantanabe-Yoshida results?} \textcolor{blue}{I think we used Verma's results in the proof of Theorem \ref{maintheoreminCMlocalringdim2} and the first example is related to  Wantanabe-Yoshida's result.}
\end{enumerate}    
\end{exmp}

%%%%%%%%%%%%%%
%%%%%%%%%%%%%%%%%%%%%%%%%%%%%%%%%%%%%

\subsection{Mixed multiplicities-like coefficients of polynomials given by length functions}
\label{subsec.p-mixedmultiplicity}

In this subsection, we will work with the following special valuation.
\begin{stp}\label{setupp-specialvaluation}
	Let $(R,\mathfrak{m},\mathbbm{k})$ be a Noetherian local domain of dimension $d$ of characteristic $p>0$ with perfect field $\mathbbm{k}$. Consider a collection of ideals $\mathcal{M}$ in $R$ such that there exists a valuation $\vartheta$ so that $\vartheta(I\cdot J)=\vartheta(I)+\vartheta(J)$ and $\vartheta(I^{[p]})=p\vartheta(I)$ for any ideals $I,J \in \mathcal{M}$. Let $\mathcal{I}(1) = \{I(1)_q\}_{q=1}^{\infty},\ldots,\mathcal{I}(1) = \{I(s)_q\}_{q=1}^{\infty}$ be $p$-families
    of ideals in such collection $\mathcal{M}$.
    
\end{stp}

\begin{exmp}
\label{ex.valofproduct}
    Let $\mathcal{M}$ be the collection of all monomial ideals and $I, J \in \mathcal{M}$. Using the notation as in Section \ref{section_about_OKvaluation}, consider any valuation $\vartheta$ given by $\textbf{a}$ whose coordinates are algebraically independent over $\ZZ$, we have that $\vartheta(f)\neq \vartheta(g)$ for any monomials $f\neq g$. We claim that $\vartheta(I\cdot J)=\vartheta(I)+\vartheta(J)$. In fact, let $a\in \vartheta(I)$, and $b\in \vartheta(J)$, then $a=\vartheta(x)$ and $b=\vartheta(y)$ for some $x\in I$ and $y\in J$. We have $a+b=\vartheta(x)+\vartheta(y)=\vartheta(xy)\in \vartheta(I\cdot J)$, hence $\vartheta(I)+\vartheta(J)\subseteq \vartheta(I\cdot J)$. Conversely, let $x=\sum\limits_{i}x_iy_i\in I\cdot J$ where $x_i\in I,y_i\in J$ are monomials, thus 
    \begin{equation*}
        \vartheta(x)=\vartheta\left(\sum\limits_{i}x_iy_i\right)=\min_{i} \{\vartheta(x_iy_i)\}=\min_{i}\{\vartheta(x_i)+\vartheta(y_i)\}\in \vartheta(I)+\vartheta(J). 
    \end{equation*}
    Also, for such valuation, since $I$ is a monomial ideal, for any $f\in I^{[p]}$, we can write $f=\sum f_ix_i^p$ where $x_i\in I$ and $f_i$ are monomials. Using the same argument as above, we get $\vartheta(f)=\vartheta(f_i)+p\vartheta(x_i) \in p\vartheta(I)$ for some $i$. This shows $\vartheta(I^{[p]})=p\vartheta(I)$.
\end{exmp}

\begin{prop}\label{prop.pbodyofproduct}
    Assume Setup \ref{setupp-specialvaluation}. For every $b \in \mathbb{N}$, and $\underline{n} = (n_1,\ldots, n_s) \in \mathbb{N}^s$, we have
    \begin{equation*}
        \Delta(S,\vartheta(\mathbf{I}(p^b)^{\underline{n}}):=\Delta(S,\vartheta(I(1)^{[p^{n_1}]}_{p^b}\cdots I(s)^{[p^{n_s}]}_{p^b}))=\sum\limits_{i=1}^{s}p^{n_i}\Delta(S,\vartheta(I(i)_{p^b}),
    \end{equation*}
    where $\Delta(S,\vartheta(I(i)_{p^b})=\bigcup\limits_{q=1}^{\infty}\left(\dfrac{1}{q}\vartheta(I(i)^{[q]}_{p^b})+\mathrm{Cone}(S)\right)$ for $i=1,\ldots,s$. 
\end{prop}

\begin{proof}
    The result can be easily obtained from the fact that 
    \begin{equation*}
        \vartheta(I(1)^{[p^{n_1}]}_{p^b}\cdots I(s)^{[p^{n_s}]}_{p^b}))=p^{n_1}\vartheta(I(1)_{p^b})+\cdots+p^{n_s}\vartheta(I(s)_{p^b}).
    \end{equation*}
\end{proof}

\begin{defn}
\label{defn.c-convexregion}
    Let $\mathcal{C}$ be a strongly convex closed $d$-dimensional cone in $\mathbb{R}^d$ with apex at the origin. A closed convex subset $\Gamma\subset \mathcal{C}$ is called $\mathcal{C}$-convex region if for any $x\in \Gamma$ and $y\in \mathcal{C}$ we have $x+y\in \Gamma$. Moreover, a convex region $\Gamma$ is cobounded if the complement $\mathcal{C}\setminus \Gamma$ is bounded. In this case the volume of $\mathcal{C}\setminus \Gamma$ is finite and is called the covolume of $\Gamma$ and denoted by $\mathrm{Covol}(\Gamma)$.
\end{defn}

\begin{lem}
    \label{lem.minkovskisum}
    Let $\mathcal{C} \subset \RR^d$ be a strongly convex closed $d$-dimensional cone with apex at the origin and suppose that there is a non-degenerate linear transformation $T$ such that $T(\mathcal{C})\subseteq \RR^d_{\ge 0}$. If $A,B$ are cobounded $\mathcal{C}$-convex regions, then so is $A+B$.
\end{lem}

\begin{proof}
    Since linear transformations preserve convexity, we can assume that $\mathcal{C} \subseteq \RR^d_{\ge 0}$. It is clear that $A+B$ is convex, cobounded, and $A+B+\mathcal{C} \subseteq A+B$. Hence, it suffices to show that $A+B$ is closed. Suppose $(c_n) \subset A+B$ is a sequence that converges to $c$. There exist sequences $(a_n)$ in $A$ and $(b_n)$ in $B$ such that $c_n=a_n+b_n$. Since all coordinates of all $b_n$ are nonnegative and $\lim c_n = c$, there exists a ball $\mathcal{B}$ in $\RR^d_{\ge 0}$ centered at $c$ such that $a_n \in \mathcal{B}$ for all $n$ large enough. Now as $\mathcal{B}$ is compact, there exists a convergent subsequence $(a_{k_n})$. Let $\lim a_{k_n}=a$, we have $a\in A$ since $A$ is closed. Moreover, $\lim b_{k_n} = \lim c_{k_n} - \lim a_{k_n} = c - a \in B$ since $B$ is also closed. Therefore, $c=a+(c-a) \in A+B$, and this implies that $A+B$ is closed as desired. 
\end{proof}

\begin{thm}
\label{thm.minsump-body}
   Assume Setup \ref{setupp-specialvaluation}. Furthermore, for every $b \in \mathbb{N}$, suppose that the $p$-body $${\Delta(S,\vartheta(I(i)_{p^b})}=\bigcup\limits_{q=1}^{\infty}\left(\dfrac{1}{q}\vartheta(I(i)^{[q]}_{p^b})+\mathrm{Cone}(S)\right)$$ is cobounded $C$-convex where $C=\mathrm{Cone}(S)$ . Then the limit $\lim\limits_{q\to \infty}\dfrac{\lambda(R/\mathbf{I}(p^b)^{q\underline{n}})}{q^d}$
   coincides with a homogeneous polynomial in $p^{n_1},\ldots,p^{n_s}$ of total degree equal to $d$. 
\end{thm}

\begin{proof}
By \cite[Corollary 5.12]{DJ18} and Proposition \ref{prop.pbodyofproduct}, we have 
\begin{align*}
     \lim\limits_{q\to \infty}\dfrac{\lambda(R/\mathbf{I}(p^b)^{q\underline{n}})}{q^d}&=[\mathbbm{k}_{\vartheta}:\mathbbm{k}] \cdot \mathrm{Covol}_{\mathbb{R}^d}(\Delta(S,\vartheta(\mathbf{I}(p^b)^{q\underline{n}})))\\
     &=[\mathbbm{k}_{\vartheta}:\mathbbm{k}] \cdot \mathrm{Covol}_{\mathbb{R}^d}(\sum\limits_{i=1}^{s}p^{n_i}\Delta(S,\vartheta(I(i)_{p^b})).
\end{align*}
Then by \cite[Theorem 2.2]{KK14}, the function $\mathrm{Covol}_{\mathbb{R}^d}(\sum\limits_{i=1}^{s}p^{n_i}\Delta(S,\vartheta(I(i)_{p^b}))$ is a homogeneous polynomial of degree $d$ in the $p^{n_i}$. It follows that $ \lim\limits_{q\to \infty}\frac{\lambda(R/\mathbf{I}(p^b)^{q\underline{n}})}{q^d}$ is also a homogeneous polynomial of degree $d$ in the $p^{n_i}$ as desired.
\end{proof}

\begin{thm}
\label{thm.minsumclosure}
   Assume Setup \ref{setupp-specialvaluation}. Suppose that there is a non-degenerate linear transformation $T$ such that $T(\mathcal{C})\subseteq \RR^d_{\ge 0}$. Furthermore, for every $b \in \mathbb{N}$, suppose that the closure of the $p$-body $\overline{{\Delta(S,\vartheta(I(i)_{p^b})}}$ is a cobounded $C$-convex region. Then the limit
   $\lim\limits_{q\to \infty}\dfrac{\lambda(R/\mathbf{I}(p^b)^{q\underline{n}})}{q^d}$
   coincides with a homogeneous polynomial in $p^{n_1},\ldots,p^{n_s}$ of total degree equal to $d$. 
\end{thm}

\begin{proof}
    First, using Lemma \ref{lem.minkovskisum} and Proposition \ref{prop.pbodyofproduct}, we get
    \[
    \overline{\Delta(S,\vartheta(\mathbf{I}(p^b)^{\underline{n}})}=\sum\limits_{i=1}^{s}p^{n_i}\overline{\Delta(S,\vartheta(I(i)_{p^b})}.
    \]
    The proof now goes along the same line with that of Theorem \ref{thm.minsump-body} once we take the closure $\overline{{\Delta(S,\vartheta(I(i)_{p^b})}}$ in the volume equalities.
\end{proof}

We are now ready to state the main theorem of this subsection.

\begin{thm}
    \label{thm.mixedMultCoeff}
    Assume Setup \ref{setupp-specialvaluation} and the conditions in either Theorem \ref{thm.minsump-body} or Theorem \ref{thm.minsumclosure}. Suppose also that $\dim (\mathrm{N}(\hat{R}))<d$. For each $b\in \NN$, set
    \[
    \mathcal{P}_{\mathbf{I}(p^b)} (p^{n_1},\ldots ,p^{n_s}) = \lim\limits_{q\to \infty}\dfrac{\lambda(R/\mathbf{I}(p^b)^{q\underline{n}})}{q^d}.
    \]
    Then there exists a homogeneous polynomial of total degree $d$ with real coefficients, denoted $\mathcal{P}_{\mathbf{I}}$, such that for all $(n_1,\ldots,n_s)\in \NN$,
    \[
    \mathcal{P}_{\mathbf{I}}(p^{n_1},\ldots ,p^{n_s}) = \lim\limits_{b\to \infty} \mathcal{P}_{\mathbf{I}(p^b)} (p^{n_1},\ldots ,p^{n_s}) = \lim\limits_{q\to \infty}\dfrac{\lambda(R/\mathbf{I}_{q\underline{n}})}{q^d}.
    \]
    
\end{thm}

In order to prove Theorem \ref{thm.mixedMultCoeff}, we are going to use similar strategies as in \cite{CSS19}. Let us write 
$$\mathcal{P}_{\mathbf{I}(p^b)} (x_1,\ldots,x_s)=\sum\limits_{i_1+\cdots+i_s=d}b_{i_1,\ldots,i_s}(p^b)x_1^{i_1}\cdots x_s^{i_s}.$$

We first show a version of \cite[Lemma 3.1]{CSS19} adapted to our context.

\begin{lem}[\cite{CSS19}, Lemma 3.1]
    \label{lem.basisOfQ^a}
    Suppose that $s,d\ge 1$, $p$ is a prime, and $a= {{s-1+d} \choose {s-1}}$. Then there exist $n_1(i),n_2(i),\ldots ,n_s(i) \in \NN$ for $1\le i \le a$ such that the set of vectors consisting of all monomials of degree $d$ in $p^{n_1(i)},p^{n_2(i)},\ldots ,p^{n_s(i)}$, $1\le i \le a$
    \[
    \left\{\left( p^{dn_1(1)},p^{(d-1)n_1(1)+n_2(1)},\ldots ,p^{dn_s(1)} \right), \ldots, \left( p^{dn_1(a)},p^{(d-1)n_1(a)+n_2(a)},\ldots ,p^{dn_s(a)} \right) \right\}
    \]
    is a $\QQ$-basis for $\QQ^a$.
\end{lem}
\begin{proof}
    The argument is similar to that of \cite[Lemma 3.1]{CSS19}. Consider the map $\Lambda: (\QQ_+)^s \rightarrow \QQ^a$ defined by $\Lambda(k_1,\ldots,k_s)=(k_1^d,k_1^{d-1}k_2,\ldots,k_s^d)$. Let $\mathcal{P} =\{p^k \ | \ k\in \NN \}$ the set of all powers of $p$. It suffices to show that the image under $\Lambda$ of the set $(\mathcal{P})^s \subset (\QQ_+)^s$ is not contained in any proper $\QQ$-linear subspace of $\QQ^a$ as it implies that we can find $a$ vectors in $(\mathcal{P})^s$ such that their images under $\Lambda$ are linearly independent, hence, form a basis for $\QQ^a$. By contradiction, suppose that $\Lambda((\mathcal{P})^s)$ is contained in a proper subspace, then there exists a nonzero linear form 
    \[
    L(y_{d,0,\ldots,0}, y_{d-1,1,0,\ldots,0}, \ldots , y_{0,0,\ldots,d})= \sum_{i_1+\ldots + i_s = d} a_{i_1,\ldots,i_s} y_{i_1,\ldots,i_s}
    \]
    on $\QQ^a$ such that $L\left(p^{dn_1},p^{(d-1)n_1+n_2},\ldots ,p^{dn_s}\right)=0$ for all $\left(p^{n_1},\ldots,p^{n_s} \right)\in (\mathcal{P})^s$. This means that the (nonzero) degree $d$ homogeneous polynomial $G(x_1,\ldots ,x_s)=L(x_1^d,x_1^{d-1}x_2,\ldots ,x_s^d)$ vanishes on $(\mathcal{P})^s$. We show that this is impossible by proving that $G \equiv 0$ using induction on $s$. Indeed, if $s=1$, since $G(x_1)$ has infinitely many roots (in $\mathcal{P}$), we must have $G(x_1)\equiv 0$. Suppose that claim is true for any homogeneous polynomial in $s-1$ variables. Write
    $G(x_1,\ldots,x_s)=\sum_{i=0}^d g_i(x_1,\ldots,x_{s-1})x_s^i$, and fix $(p^{n_1},\ldots ,p^{n_{s-1}}) \in (\mathcal{P})^{s-1}$. The polynomial $G(p^{n_1},\ldots ,p^{n_{s-1}},x_s)\in \QQ[x_s]$ has infinitely many solutions, hence, $G(p^{n_1},\ldots ,p^{n_{s-1}},x_s) \equiv 0$, therefore, for all $0\le i \le d$, $g_i(p^{n_1},\ldots ,p^{n_{s-1}})=0$. Since $(p^{n_1},\ldots ,p^{n_{s-1}})$ is arbitrary chosen in $(\mathcal{P})^{s-1}$, by induction, $g_i\equiv 0$ for all $0\le i\le d$, thus, $G\equiv 0$ as desired.
\end{proof}

\begin{prop}[\cite{CSS19}, Corollary 4.4]
\label{prop.limitOfCoeff}
    For all $i_1,\ldots,i_s\in \mathbb{N}$ with $i_1+\cdots+i_s=d$, the limit
        $\displaystyle b_{i_1,\ldots,i_s}=\lim_{b\to \infty}b_{i_1,\ldots,i_s}(p^b)$
    exists in $\RR$.
\end{prop}

\begin{proof}
    The proof is the same as in \cite[Corollary 4.4]{CSS19} since we can use Lemma \ref{lem.basisOfQ^a} to show \cite[Lemma 3.2]{CSS19} and use Theorem \ref{maintool} in place of \cite[Proposition 4.3]{CSS19}.
\end{proof}

We are now ready to prove Theorem \ref{thm.mixedMultCoeff}.

\begin{proof}[Proof of Theorem \ref{thm.mixedMultCoeff}]
    We define a homogeneous polynomial 
\begin{equation*}
    \mathcal{P}_{\mathbf{I}}(x_1,\ldots,x_s)=\sum\limits_{i_1+\cdots+i_r=d}b_{i_1,\ldots,i_r}x_1^{i_1}\cdots x_s^{i_s}.
\end{equation*}
 For any $(n_1,\ldots,n_s)\in \NN$, by Theorem \ref{thm.maintool}, we have
    \[
    \lim_{b\to \infty} \mathcal{P}_{\mathbf{I}(p^b)} (p^{n_1},\ldots ,p^{n_s}) = \lim_{b\to \infty} \lim_{q\to \infty}\dfrac{\lambda(R/\mathbf{I}(p^b)^{q\underline{n}})}{q^d} = \lim_{q\to \infty}\dfrac{\lambda(R/\mathbf{I}_{q\underline{n}})}{q^d},
    \]
    Moreover, by Proposition \ref{prop.limitOfCoeff}, $\displaystyle \lim_{b\to \infty} \mathcal{P}_{\mathbf{I}(p^b)} (p^{n_1},\ldots ,p^{n_s}) = \mathcal{P}_{\mathbf{I}}(p^{n_1},\ldots ,p^{n_s})$ and this concludes the proof.
\end{proof}

\end{document}